\documentclass[a4paper,11pt,reqno]{amsart}

\usepackage{scalerel,stackengine}
\stackMath
\usepackage{amsmath,amssymb,amsthm,mathrsfs,amsfonts,dsfont,stmaryrd} 

\usepackage{multirow}
\usepackage{amscd}   
\usepackage{fullpage}
\usepackage[all]{xy}  
\usepackage{graphicx}
\usepackage{mathrsfs,enumitem}
\usepackage[titletoc,toc,title]{appendix}
\usepackage{xcolor}  
\usepackage[normalem]{ulem}
\usepackage{diagbox}
\usepackage{caption}

\usepackage[colorlinks = true,
            linkcolor = blue,
            urlcolor  = blue,
            citecolor = blue,
            anchorcolor = blue]{hyperref}
\usepackage{graphicx}
 \usepackage{color}
 \usepackage{tikz, framed}
\usetikzlibrary{arrows,shapes,matrix,decorations.pathmorphing}
\usetikzlibrary{backgrounds}
\usepackage[english]{babel} 
\newtheorem{Theorem}{Theorem} 

\newtheorem{Lemma}[Theorem]{Lemma}
\newtheorem{Proposition}[Theorem]{Proposition}
\newtheorem{Corollary}[Theorem]{Corollary}

\theoremstyle{definition}
\newtheorem{Definition}[Theorem]{Definition}
\newtheorem{Remark}[Theorem]{Remark}
\newtheorem{Example}[Theorem]{Example}
\newtheorem{Notation}[Theorem]{Notation}

 \newcommand{\kB}{\mathcal{B}}
 \newcommand{\kL}{\mathcal{L}}
 
  \newcommand{\mL}{\mathcal{L}}
    \newcommand{\mC}{\mathcal{C}}

 \newcommand{\NN}{\mathbb{N}}
 \newcommand{\ZZ}{\mathbb{Z}}

 \newcommand{\FF}{\mathbb{F}}
 \newcommand{\mP}{\mathcal{P}}
  \newcommand{\D}{\mathcal{D}}
 \newcommand{\mB}{\mathcal{B}}
  
 \newcommand{\Fq}{\mathbb{F}_q}

 \newcommand{\Fqm}{{\mathbb{F}_{q^m}}}

\newcommand\qbin[3]{\left[\begin{matrix} #1 \\ #2 \end{matrix} \right]_{#3}} 

\DeclareMathOperator{\rk}{rank}
\DeclareMathOperator{\cl}{cl}
\DeclareMathOperator{\col}{colsp}

\begin{document}  

\title{Weighted Subspace Designs from $q$-Polymatroids}

\author[Eimear Byrne]{Eimear Byrne}
\address{School of Mathematics and Statistics, University College Dublin, Belfield, Ireland}
\curraddr{}
\email{ebyrne@ucd.ie}
\thanks{}

\author[Michela Ceria]{Michela Ceria}
\address{
Dept. of Mechanics, Mathematics \& Management, Politecnico di Bari, Italy
}
\curraddr{  Via Orabona 4 - 70125 Bari - Italy}
\email{michela.ceria@gmail.com}
\thanks{}

\author[Sorina Ionica]{Sorina Ionica}
\address{
University of Picardie Jules Verne 
}
\curraddr{33 rue Saint Leu Amiens, 80039, France}
\email{sorina.ionica@u-picardie.fr}
\thanks{}

\author[Relinde Jurrius]{Relinde Jurrius}
\address{Faculty of Military Sciences, Netherlands Defence Academy, The Netherlands}
\curraddr{}
\email{rpmj.jurrius@mindef.nl}
\thanks{}

\begin{abstract}
The Assmus-Mattson Theorem gives a way to identify block designs arising from codes. This result was broadened to matroids and weighted designs by Britz et al. in 2009. In this work we present a further two-fold generalisation: first from matroids to polymatroids and also from sets to vector spaces. To achieve this, we study the characteristic polynomial of a $q$-polymatroid and outline several of its properties. We also derive a MacWilliams duality result and apply this to establish criteria on the weight enumerator of a $q$-polymatroid for which dependent spaces of the $q$-polymatroid form the blocks of a weighted subspace design.
\end{abstract}



\maketitle
\section{Introduction}

The characteristic polynomial of a matroid is a well studied object. It was first introduced as a matroid generalisation of the chromatic polynomial of a graph. It arises in critical problems, analyses of the Tutte polynomial, and is the subject of numerous identities \cite{kung1996critical} 
For a thorough treatment of the subject see \cite{welsh, zaslavsky1987mobius}, for example. 

In combinatorics, the concept of a $q$-analogue can be viewed as a generalisation from sets to vector spaces. Recently, the $q$-analogue of a matroid has been studied \cite{JP18}. A generalisation of this is a $q$-polymatroid \cite{GLJ,gorla2019rank,shiromoto19}.

Similar to classical matroids, there are many interesting connections between $q$-(poly)matroids and rank-metric codes. In this paper we develop the theory of the characteristic polynomial of a $q$-polymatroid. We show the relation between the characteristic polynomial of a polymatroid and its dual, establishing a MacWilliams-like identity for $q$-polymatroids. In a similar line of research, Shiromoto \cite{shiromoto19} established a $q$-analogue of Greene's theorem. 

Another motivation to study the characteristic polynomial is to establish a $q$-analogue of the Assmus-Mattson Theorem \cite{AM69}. This theorem gives a criterion for identifying a $t$-design as a collection of supports of codewords of fixed weight in a linear code.
Since its publication in 1969 it has seen a number of generalisations \cite{byrne2019assmus,MT18} and has been used widely to obtain new constructions of designs \cite{DL17,SKH04}. In one of these results \cite{BRS2009}, the authors define a weighted $t$-design as a generalisation of a classical $t$-design and give criteria for identifying such an object among the dependent sets of a matroid of a fixed cardinality. A weighted $t$-design is a collection of subsets $\mB$ of a fixed cardinality $k$ chosen from an $n$-set of points $\mP$ together with a function $f$ defined on $\mB$ such that for any $t$-set $T \subset \mP$ the sum $\sum_{B \in {\mB\,:\,T}\subset B} f(B)$ is independent of $T$.
In the case that $f(B)=1$ for every block $B \in \mB$,  
the weighted $t$-design is an ordinary design.

In this paper, we generalise the results of \cite{BRS2009} to $q$-polymatroids, which is a two-fold generalisation: first from matroids to polymatroids and also from sets to vector spaces. Hence the results presented here give a $q$-analogue of their result. The $q$-analogue of a weighted $t$-design is a weighted subspace design; in the definition shown above we replace the collection of subsets $\mB$ with a collection of subspaces of a fixed dimension $k$ and $T$ with a $t$-dimensional subspace.

In Section \ref{sec:qpoly} we study $q$-polymatroids and the properties of $q$-polymatroids that are necessary for this work.
In Section \ref{sec:charpoly} we outline properties of the characteristic polynomial of a $q$-polymatroid that will be used later and in Section \ref{sec:codes} look at the case of $q$-polymatroids arising from matrix codes. In Section \ref{sec:macwill} we give a version of the MacWilliams duality result for $q$-polymatroids. In Section \ref{sec:designs} we give criteria for identifying when the dependent spaces of a $q$-polymatroid are the blocks of a weighted $t$-subspace design. \\

\begin{Notation}
Throughout, we let $n$ denote a fixed positive integer and we will let $q$ denote a fixed prime power.
 We let $E$ denote an $n$-dimensional vector space over the finite field $\mathbb{F}_q$ of order $q$. We let $\mathcal{L}(E)$ denote the lattice of all subspaces of $E$, ordered  by inclusion, which we denote by $\leq $. We will write $U< V$ for $U,V\leq E$ if $U$ is strictly contained in $V$. The join of a pair of subspaces is their vector space sum and the meet of a pair of subspaces is their intersection.
\end{Notation}

\section{{\em q}-Polymatroids}\label{sec:qpoly}

$q$-Polymatroids and their connections to linear codes were introduced in \cite{gorla2019rank} and \cite{shiromoto19}. Their properties have been further developed in \cite{GLJ}. In our presentation, we will not assume that $q$-polymatroids are {\em representable}, that is, we will not assume that the $q$-polymatroids under consideration here are constructed from rank-metric codes over $\FF_q$.  
We use the following definition of a $q$-polymatroid from \cite{shiromoto19}, since it suits our purposes to have an integer valued function in what follows. 

\begin{Definition}
A $(q,r)$-polymatroid is a pair $M=(E, \rho)$ for which $r \in \NN_0$ and $\rho$ is a function $\rho: \kL(E) \longrightarrow {\mathbb N}_0$ satisfying the following axioms. 
\begin{itemize}
	\item[(R1)] For all $A \leq E$, $0\leq \rho(A) \leq r\dim A$.
	\item[(R2)] For all $A,B \leq E$, if $A\leq B${, then} $\rho(A)\leq \rho(B)$.
	\item[(R3)] For all $A,B \leq E$, $\rho(A+B)+\rho(A\cap B)\leq \rho(A)+\rho(B)$.
\end{itemize}
\end{Definition}
Every $(q,r)$-polymatroid is also a $(q,r')$-polymatroid for any $r'\geq r$. Hence, all  the definitions below involving $r$ depend on the choice of $r$.
If it is not necessary to specify $r$, we will simply refer to such an object as a $q$-polymatroid. If we need to specify the $q$-polymatroid $M$, we denote its rank function by $\rho_M$. Note that a $(q,1)$-polymatroid is a $q$-matroid. Conversely, if we consider a $q$-matroid as a $(q,r)$-polymatroid, we will always take $r=1$.
In order to stress in a stronger way the distinction between matroids and their $q$-analogues, we use the terminology ``classical matroids'' for matroids.

Recall that a lattice isomorphism between a pair of lattices $(\kL_1,\vee_1,\wedge_1),(\kL_2,\vee_2,\wedge_2)$ is a bijective function $\varphi:\kL_1\longrightarrow\kL_2$ that preserves the meet and join, that is, for all $x,y\in\kL_1$ we have that $\varphi(x\wedge_1 y)=\varphi(x)\wedge_2\varphi(y)$ and $\varphi(x\vee_1 y)=\varphi(x)\vee_2\varphi(y)$.  {
It is well known that reversing the ordering of a lattice gives again a lattice, with the meet and join interchanged. Combining this operation with a lattice isomorphism gives a lattice anti-isomorphism. Formally, a lattice anti-isomorphism between a pair of lattices is a bijective function $\psi:\kL_1\longrightarrow\kL_2$ that is order-reversing and interchanges the meet and join, that is, for all $x,y\in\kL_1$ we have that $\psi(x\wedge_1 y)=\psi(x)\vee_2\psi(y)$ and $\psi(x\vee_1 y)=\psi(x)\wedge_2\psi(y)$. See \cite[Pages 3--4]{Birkhoff}.}
We hence define a notion of equivalence {and duality} between $q$-polymatroids.

\begin{Definition}\label{def:lat-equiv}
Let $E_1,E_2$ be $\FF_q$-linear spaces. Let $M_1=(E_1,\rho_1)$ and $M_2=(E_2,\rho_2)$ be $q$-polymatroids. We say that $M_1$ and $M_2$ are \textbf{lattice-equivalent} if there exists a lattice isomorphism $\varphi:\kL(E_1)\longrightarrow \kL(E_2)$ such that $\rho_1(A)=\rho_2(\varphi(A))$ for all $A\leq E_1$. In this case we write $M_1 \cong M_2$.
\end{Definition}

{
\begin{Notation}
\label{def:bform}
Let $F$ be an $\mathbb{F}_q$-vector space. We denote by $\perp\!\!(F)$ a fixed anti-isomorphism on $\kL(F)$, which we require to be an involution.  For any subspace $U\leq F$ we denote by $U^{\perp(F)}$ the image of $U$ under $\perp\!\!(F)$, which we call the \textbf{dual} of $U$ in $F$.  Note that since an anti-isomorphism preserves the length of intervals, we have for any $U\leq F$ that $\dim(U^{\perp(F)})=\dim(F)-\dim(U)$.
 In the case $F=E$, we drop the $(E)$ and we simply write  $\perp:=\perp\!\!(E)$. For any subspace $U\leq E$, we write $U^\perp:=U^{\perp(E)}$.

\end{Notation}
}

\begin{Remark}\label{rem:LatEqVsEq}
   {The notion of lattice-equivalence in Definition~\ref{def:lat-equiv} is not the same as the definition of equivalence of $q$-polymatroids given in \cite{GLJ} and \cite{gorla2019rank}. Indeed, in ~\cite{GLJ} and~\cite{gorla2019rank} two $q$-polymatroids $(E_1,\rho_1)$ and $(E_2,\rho_2)$ are said to be {\bf equivalent} if there exists an $\FF_q$-linear isomorphism $\tau: E_1 \longrightarrow E_2$ such that $\rho_1(A)=\rho_2(\tau(A))$ for all $A \leq E_1$.} Since every vector space isomorphism induces a lattice isomorphism, equivalence implies lattice-equivalence for $q$-polymatroids. {In particular, every non-degenerate symmetric bilinear form $b_F$ on $F$ induces a lattice anti-isomorphism, hence our definition of dual implies the usual definition of orthogonal complement for $q$-polymatroids.}
\end{Remark}

The dual $q$-polymatroid was defined in \cite{gorla2019rank,shiromoto19}.

\begin{Definition}\label{dualit}
Let $M=(E,\rho)$ be a $(q,r)$-polymatroid. For every subspace $A\leq E$, define $\rho^*(A):=r\dim(A)-\rho(E)+\rho(A^{\bot})$.
Then $M^*:=(E,\rho^*)$ is a $(q,r)$-polymatroid called the {{\bf lattice-dual}} of $M$.
\end{Definition}

We call $M^*$ the \textbf{dual} of $M$.
As noted in \cite{GLJ}, the definition of the dual of $M$ depends on the choice of
{anti-isomorphism on $E$}, but all such choices yield equivalent duals.
{We prove this for our more general notions of lattice-equivalence and lattice-duality. The following is a generalisation of \cite[Theorem 2.8]{GLJ}.
\begin{Lemma}\label{DualEquiv}
Let $M=(E,\rho)$ be a $(q,r)$-polymatroid and let $M'=(E,\rho')$ be a $(q,r)$-polymatroid that is lattice-equivalent to $M$. Let $\perp,\hat{\perp}$ be a pair of anti-isomorphisms on $\kL(E)$. 
Let $M^*$ and $M^{\hat{*}}$ be the duals of $M$ with respect to $\perp$ and $\hat{\perp}$, respectively. Then $M^*\cong M'^*$ and $M^*\cong M^{\hat{*}}$.
\end{Lemma}
\begin{proof}
For the first part, notice that the proof of \cite[Proposition 3.7]{gorla2019rank} carries over directly from $\mathbb{F}_q$-isomorphisms to lattice-(anti-)isomorphisms. We include it here for completeness. Let $\varphi:\kL(E)\longrightarrow\kL(E)$ be the isomorphism such that $\rho(A)=\rho'(\varphi(A))$ for all $A\subseteq E$. Let $\psi:\kL(E)\longrightarrow\kL(E)$ be the isomorphism 
$\psi:=\,\perp\circ\ \varphi\ \circ\perp$, so $\varphi\ \circ\perp\,=\,\perp\circ\ \psi$.
Then
\begin{align*}
\rho^*(A) & =r\dim(A)-\rho(E)+\rho(A^\perp) \\
 & = r\dim(A)-\rho'(E)+\rho'(\varphi(A^\perp)) \\
 & = r\dim(\psi(A))-\rho'(\psi(E))+\rho'(\psi(A)^\perp) \\
 & = \rho'^*(\psi(A)).
\end{align*}
This shows that $M^*\cong M'^*$. For the second statement we proceed in a similar way.  
Let $\phi:\kL(E)\longrightarrow\kL(E)$ be the lattice isomorphism $\phi:=\perp \circ \hat{\perp}$, so $\perp=\phi\circ \hat{\perp}$. Then
\begin{align*}
\rho^*(A) & =r\dim(A)-\rho(E)+\rho(A^\perp) \\
& = r\dim(\phi(A))-\rho(\phi(E))+\rho(\phi(A^{\hat{\perp}})) \\
& = \rho^{\hat{*}}(\phi(A)).
\end{align*}
This shows that $M^*\cong M^{\hat{*}}$.
\end{proof}
Note that $M^{**}=M$ is an equality, because we assume the anti-isomorphism $\perp$ to be an involution.}

It is easy to see that for a map $\rho: \kL(E) \longrightarrow {\mathbb N}_0$ satisfying the axioms (R1)-(R3), the restriction of that map to $\kL(T)$, for each subspace $T \leq E$,  also yields a $q$-polymatroid.  

\begin{Definition}
Let $M=(E,\rho)$ be a $(q,r)$-polymatroid and let $T\leq E$. 
For every subspace $A\leq T$, define $\rho_{M|T}(A):=\rho(A)$.
Then $M|T:=(T,\rho_{M|T})$ is a $(q,r)$-polymatroid called the {\bf restriction} of $M$ to $T$.
\end{Definition}

{Another way to construct a new $q$-polymatroid from an existing one is via contraction. It was proven in \cite[Theorem 5.2]{GLJ} that this gives in fact a $q$-polymatroid. \\}

\begin{Definition}
Let $M=(E, \rho)$ be a $(q,r)$-polymatroid and let $T\leq E$.
We define the map
\[ \rho_{M/T}: \kL(E/T) \longrightarrow \ZZ \]
 via $\rho_{M/T}(A/T) = \rho(A) - \rho(T)$.
Then $M/T:=(E/T,\rho_{M/T})$ is a $(q,r)$-polymatroid called the   {\bf contraction} of $T$ {\bf from} $M$. 
\end{Definition}

It will sometimes be more convenient for us to use the slightly less commonly used definition of contraction \emph{to} a subspace.

\begin{Definition}\label{Cont}
    Let $M=(M,\rho)$ be a $(q,r)$-polymatroid and 
	let $X \leq E$. We denote by $M.X$ the $q$-polymatroid 
	$M.X:=(E / X^\perp, \rho_{E / X^\perp})$. We call
	$M.X$ the {\bf contraction} of $M$ {\bf to} $X$.
\end{Definition}

In the language of classical matroids, the contraction of $M$ to $X$ is the contraction
of $E-X$ from $M$, that is $M.X = M/(E-X)$ (see \cite[Chapter 3]{oxley}). In the $q$-analogue we have 
$M.X:=M/X^\perp$.

The following duality result is a straightforward extension of \cite[Theorem 60]{JP18}. 
It relates the contraction of a subspace from a $q$-polymatroid to a restriction of its dual $q$-polymatroid. We will make good use of this in Section \ref{sec:designs}, where we give a construction of weighted subspace designs from $q$-polymatroids. 		

\begin{Lemma}\label{lem:dualMTperp}
	Let $M=(E,\rho)$ be a $(q,r)$-polymatroid and let $T$ be a subspace of $E$. Then,
	$$M^*/T \cong ({M|T^\perp})^* \text{ and } (M/T)^* \cong M^*|T^\perp. $$
\end{Lemma}
\begin{proof}	
  Let $\phi: {\mathcal L}(E/T) \longrightarrow {\mathcal L}(T^\perp)$ be defined by $\phi(X/T)= (X^\perp)^{\perp(T^\perp)}$, for each $X \leq E$ such that $T \leq X$ (in which case $X^\perp \leq T^\perp$).
  This map is the composition of two anti-isomorphisms: the anti-isomorphism between intervals  $[T,E]$ and $[0,T^\perp]$ induced by $\perp\!\!(E)$, followed by the anti-isomorphism $\perp\!\!(T^\perp)$ on $\kL(T^\perp)$. Hence $\phi$ is a lattice isomorphism.

 Let $A$ be a subspace of $E$ satisfying $T \leq A \leq E$. 
  We claim that $\rho_{M^*/T}(A/T) = (\rho_{{M|T^\perp}})^*(\phi(A/T))$. Indeed, we have that:
  \begin{eqnarray*}
  	  \rho_{M^*/T}(A/T)& = & \rho^*(A)-\rho^*(T) \\
  	                & = & r\dim(A)-\rho(T^\perp)+\rho(A^\perp) -r\dim(T) \\
  	                & = & r\dim(A/T)-\rho_{{M|T^\perp}}(T^\perp)+\rho_{{M|T^\perp}}(A^\perp)\\
  	                & = & r\dim(\phi(A/T)) - \rho_{{M|T^\perp}}(T^\perp) + \rho_{{M|T^\perp}}(\phi(A/T)^{\perp(T^\perp)})\\
  	                & = & (\rho_{{M|T^\perp}})^*(\phi(A/T)).
  \end{eqnarray*}
  This shows that $M^*/T \cong ({M|T^\perp})^*$. That $(M/T)^* \cong M^*|T^\perp$ holds can be seen by replacing $M$ with $M^*$ in the previous identity, taking duals {and applying Lemma \ref{DualEquiv}}.
\end{proof}

\begin{Remark}
   In fact, the above result holds even in terms of equivalence in the stronger sense 
   \cite[Definition 2.6 (b)]{GLJ}, and not only lattice-equivalence, as was shown in Theorem 5.3 
   of the same paper.
   Note that in establishing the equivalence of these $q$-polymatroids, the vector space isomorphism 
   depends on the choice of the bilinear form arising in the construction of the lattice isomorphism.    
\end{Remark}

Having established duality, restriction and contraction in terms of the rank function, we now introduce independent spaces.

\begin{Definition}
Let $I\leq E$ and let $M=(E,\rho)$ be a $(q,r)$-polymatroid. We say that $I$ is an {\bf independent space} of $M$ if $\rho(I)=r\dim I$. A subspace that is not independent is called a {\bf dependent space} of $M$. We call $C\leq E$ a {\bf circuit} of $M$ if it is a minimal dependent space with respect to inclusion. We call $T \leq E$ a {\bf cocircuit} of $M$ if it is a circuit of $M^*$.
\end{Definition}

For $q$-matroids, the following result is (I2) of the independence axioms (see \cite[Definition 7]{bcj}). We show that this holds for $q$-polymatroids.

\begin{Lemma}\label{lem:I2}
Let $M=(E,\rho)$ be a $(q,r)$-polymatroid and let $I \leq E$ be an independent space of $M$. Then every subspace of $I$ is independent.
\end{Lemma}
\begin{proof}
 Since $I$ is independent, we have $\rho(I)=r\dim(I)$. Let $J,J'$ be subspaces of $I$ such that $I$ is a direct sum of $J$ and $J'$. 
 By (R1) and applying semimodularity (R3) to $J$ and $J'$ we get
\begin{align*}
r\dim(J) + r\dim(J')&\geq \rho(J) + \rho(J') \geq \rho(J+J')+\rho(J\cap J')\\
&= \rho(I) = r\dim(I) = r(\dim(J) + \dim(J')).
\end{align*}
Since $\rho(J)\leq r\dim(J)$ and $\rho(J') \leq r \dim(J')$ we must have that $\rho(J)=r\dim(J)$ and $\rho(J')=r\dim(J')$ and the result follows.
\end{proof}

From the above lemma, it follows that  $C\leq E$ is a circuit of a $q$-polymatroid if it is a dependent space whose proper subspaces are all independent. If $J$ is a maximal independent subspace of $A \leq E$, then $\rho(A)=\rho(J)$; for $q$-matroids this follows directly from the cryptomorphism between independent spaces and rank \cite[Theorem 8]{JP18}, for $q$-polymatroids this was shown in \cite[Theorem 4.2]{GLJindep}.

{The next lemma considers what happens to independent spaces and circuits under contraction of an independent space.}

\begin{Lemma}\label{lem:MtoM/T}
    Let $M=(E,\rho)$ be a $(q,r)$-polymatroid and let $I \leq E$ be an independent space of $M$.
	Let $I \leq A \leq E$.
	Then $A$ is independent in $M$ if and only if $A/I$ is independent in $M/I$.
	Moreover, if $A$ is a circuit in $M${, then} $A/I$ is a circuit in ${M/I}$.
\end{Lemma}
\begin{proof}
	Let $A$ be independent in $M$. Then
	$$r\dim(A/I) = r\dim(A)-r\dim(I)=\rho(A)-\rho(I) = \rho_{M/I}(A/I),$$
	hence $A/I$ is an independent space of $M/I$.
	Conversely, if $A/I$ is independent in $M/I${, then}
	$$r\dim(A)-r\dim(I)=r\dim(A/I)=\rho_{M/I}(A/I) = \rho(A)-\rho(I) =\rho(A)-r\dim(I),$$
	so $\rho(A)=r\dim(A)$.
	
	Let $A$ be a circuit in $M$. Any proper subspace of $A/I$ has the form $B/I$ for some
	unique $I \leq B < A$. Since $A$ is a circuit, $A/I$ is a dependent space
	in $M/I$, and $B$ is an independent space of $M$. 
	Therefore $B/I$ is independent and so $A/I$ is a circuit
	of $M/I$.
\end{proof}

We conclude this section with some examples. These are from \cite[Example 4]{JP18} and \cite[Example 4.8(a)]{GLJ}.

\begin{Example}
[The uniform $q$-Matroid] \label{exunif}
	Let $k$ be a positive integer, $k \leq n$. 
	The  {\em uniform $q$-matroid}  is the $q$-matroid $M=(E,\rho)$  with rank function defined as follows:
	\[
	  \rho(U):= \left\{ \begin{array}{cl}
	  	  \dim(U) & \text{ if } \dim(U) \leq k, \\
	  	  k       & \text{ if } \dim(U) > k.
	  \end{array}\right.
	\] 
	We denote this $q$-matroid by $U_{k,n}$.
\end{Example}

\begin{Example}
[The V\'{a}mos $q$-Matroid ]\label{exVamos}
This $q$-matroid is constructed over $\kL(\FF_q^8)$. Choose the canonical basis for $\FF_q^8$ denoted by $e_1,\ldots,e_8$. Consider the following collection of subspaces. 
\[ \mC:=\{  \langle e_1,e_2,e_3,e_4\rangle,
\langle e_1,e_2,e_5,e_6\rangle,
\langle e_3,e_4,e_5,e_6\rangle,
\langle e_3,e_4,e_7,e_8\rangle,
\langle e_5,e_6,e_7,e_8\rangle\}.
\]
For each $A\leq \FF_q^8$, we define $\rho(A)$ as follows:

\[
\rho(A):=\left\{\begin{array}{cl}
    \dim(A) & \textrm{ if } \dim(A)\leq 3,  \\
    3 & \textrm{ if } A \in \mC,\\
4 &  \textrm{ if } \dim(A)\geq 4 \textrm{ and } A \notin \mC.
\end{array}\right.
\]
It can be shown that $\rho$ is the rank function of a $q$-matroid whose set of circuits of minimum dimension is the set $\mC$.
\end{Example}

\section{Characteristic Polynomial of a {\em q}-Polymatroid}\label{sec:charpoly}

In this section, we introduce the characteristic polynomial of a $q$-polymatroid. This polynomial and its properties are well-studied in the case of a classical polymatroid
\cite{kung1996critical,whittle}
, in which case its coefficients are the M\"{o}bius values of the lattice of subsets of $\{1,\ldots,n\}$. In the $q$-polymatroid case the underlying lattice is the lattice of subspaces of $E$. We will use the characteristic polynomial to obtain a version of the MacWilliams identities for $q$-polymatroids.

\subsection{The M\"obius Function of a Lattice}

Throughout this paper we will use the M\"obius function  (see, e.g., ~\cite[Chapter 25]{Lint}), which is fundamental to the definition of a characteristic polynomial. We recall some basic results.

Let $(P,\leq)$ be a partially ordered set. The M\"obius function for $P$ is defined via the recursive formula
\begin{align*}
\mu(x,y) := \left\{ 
        \begin{array}{cl}
        1 & \text{ if } x=y,\\
        -\sum_{x\leq z<y}\mu(x,z) & \text{ if } x<y, \\
        0 & \text{otherwise}.
        \end{array}
          \right.
\end{align*}
\begin{Lemma}[M\"obius Inversion Formula]\label{inversionformula}
  Let $f,g,h:P\longrightarrow \mathbb{Z}$ be any 3 functions on a poset $P$. Then, we have
	\begin{enumerate}
		\item $\displaystyle f(x)=\sum_{x\leq y}g(y)$ for all $x\in P$ if and only if $\displaystyle g(x)=\sum_{x\leq y}\mu(x,y)f(y)$ for all $x \in P$,
		\item  $\displaystyle f(x)=\sum_{x\geq y}h(y)$  
		for all $x\in P$
		if and only if $\displaystyle h(x)=\sum_{x\geq y}\mu(y,x)f(y)$ for all $x\in P$.
	\end{enumerate} 
\end{Lemma}

For the subspace lattice of $E$ and for two subspaces $U$ and $V$ of dimensions $u$ and $v$, we have that
\begin{equation*}\label{eq:mobinv}
\mu\left(U,V\right)=\left\lbrace
\begin{array}{cl}
(-1)^{v-u}q^{\binom{v-u}{2}}&\mbox{ if }U\leq V\\
\\
0 & \mbox{ otherwise.}
\end{array}
\right.
\end{equation*}

\begin{Definition}
    Given a pair of nonnegative integers $a$ and $b$, the $q$-\textbf{binomial} or \textbf{Gaussian coefficient} counts the number of $b$-dimensional subspaces of an $a$-dimensional subspace over $\FF_q$ and is given by:
$$\qbin{a}{b}{q}:=\prod_{i=0}^{b-1}\frac{q^a-q^i}{q^b-q^i},$$
if $a \geq b$ and is zero if $a<b$.
\end{Definition}

We will use the following identities 
\begin{equation}\label{eq:qbinid}
   \qbin{a}{b}{q} = \qbin{a}{a-b}{q} \text{ and } \qbin{a}{b}{q}\qbin{b}{c}{q} = \qbin{a}{c}{q} \qbin{a-c}{a-b}{q}.
\end{equation}
See also \cite{qdist}, for example, for a comprehensive account of the properties of Gaussian coefficients.

\begin{Lemma}\label{counting_subspace}
	Let $I, J$ be subspaces of $E$ of dimensions $i$ and $j$, respectively, satisfying  $I \cap J=\{0\}$ and $i+j \leq k$. Then, the number of $k$-dimensional subspaces of $E$ that contain $I$ and meet trivially with $J$ is
    \[ \sum_{s=0}^{j}(-1)^{s}q^{\binom{s}{2}} \qbin{j}{s}{q}\qbin{n-i-s}{k-i-s}{q}=	q^{j(k-i)}\qbin{n-i-j}{k-i}{q} \]
where $n$ is the dimension of E.	
\end{Lemma}
In the above Lemma \ref{counting_subspace}, we note another identity. 
That the right-hand side counts the number of $k$-dimensional subspaces of $E$ that contain $I$ and meet trivially with $J$ was already observed, for example, in \cite[Lemma 3]{del76}, but is generally well-known. 
That this number is also given by the left-hand side formula can be established using M\"obius inversion.

\subsection{The Characteristic Polynomial}

 We now introduce the characteristic polynomial. First, we give another definition, which originates in weight enumeration of linear codes.

\begin{Definition}  \label{funzionel}
Let $M$ be a $(q,r)$-polymatroid with ground-space $E$. For each $A \leq E$ we define 
\[
    \ell_M(A):=\rho_M(E)-\rho_M(A).
\]
\end{Definition}

By the definition of the rank function of a $q$-polymatroid, for each subspace $A$ of $E$ 
we see that $\ell_M(A)$ is non-negative integer in $\{0,\ldots,\rho_M(E)\}$.

{
\begin{Notation}
    For the remainder, we fix $r$ to be a positive integer and we let $M$ denote a fixed $(q,r)$-polymatroid $M=(E,\rho)$. We write $\ell:=\ell_M$ and $\rho:=\rho_M$.
    For the dual $q$-polymatroid, we write $\ell^*:=\ell_{M^*}$ and $\rho^*:=\rho_{M^*}$.
    \end{Notation}}

\begin{Definition}
	The {\bf characteristic polynomial} of $M$ is the polynomial in $\ZZ[z]$ defined by
	$$ p(M;z) := \sum_{{X\,:\,X }\leq E} \mu(0,X) z^{\ell(X)},$$
	where $\mu$ is the M\"obius function of the subspace lattice.
\end{Definition}

 For the case $E=\{0\}$, we have $p(M;z)=1$. If $E \neq \{0\}$, then $p(M;1)=0$ and so, unless $p(M;z)$ is identically zero, $z-1$ is a factor of it in $\ZZ[z]$.

For a $(q,r)$-polymatroid $M$, we have $$p(M;z) := \sum_{j=0}^n (-1)^j q^{\binom{j}{2}} \sum_{X\leq E, \dim(X)=j} z^{\ell(X)}. $$ 

\begin{Example}
We calculate the characteristic polynomial of the V\'{a}mos $q$-matroid of Example \ref{exVamos}. From the rank function it follows that:
\[
\ell(X)=\left\{\begin{array}{cl}
    4-\dim(X) & \textrm{ if } \dim(X)\leq 3,  \\
    1 & \textrm{ if } X \in \mC,\\
0 &  \textrm{ if } \dim(X)\geq 4 \textrm{ and } X \notin \mC.
\end{array}\right.
\]
 We treat the calculations of the coefficients by the powers of $z$. For the coefficient of $z^4$ we only have $X\leq E$ with $\dim X=0$,  i.e., the zero space. Then $\mu(0,X)=\mu(0,0)=1$ and we get the term $z^4$.
For $z^3$ and $z^2$ we get 
\[ \sum_{\dim X=1} \mu(0,X)z^{\ell(X)}=-\qbin{8}{1}{q}z^3,\qquad \sum_{\dim X=2} \mu(0,X)z^{\ell(X)}=q\qbin{8}{2}{q}z^2. \]
For the coefficient of $z$ we have to consider the five circuits of dimension $4$ and all spaces of dimension $3$, which is
\[ 5q^6-q^3\qbin{8}{3}{q}. \]
Finally, the constant term is determined by all spaces of dimension $4$ that are not circuits, plus all spaces of higher dimension:
\begin{small}
\begin{align*}
    q^6\left(\qbin{8}{4}{q}-5\right)-q^{10}\qbin{8}{5}{q}+q^{15}\qbin{8}{6}{q}-q^{21}\qbin{8}{7}{q}+q^{28}
    =q^6\left( \qbin{8}{4}{q}-5-q^{4}\qbin{8}{3}{q}+q^{9}\qbin{8}{2}{q}-q^{15}\qbin{8}{1}{q}+q^{22} \right).
\end{align*} 
\end{small}
Adding all terms gives the characteristic polynomial of the V\'{a}mos $q$-matroid.
For example, for $q=2$, we have $p(M;z)=z^4-255z^3+21590z^2-776920z+755584 = (z-1)(z^3-254z^2+21336z-755584)$.
\end{Example}

It is easily checked that the characteristic polynomial is an invariant of the lattice-equivalence class of a matroid.

\begin{Lemma}
    Let $E_1,E_2$ be $\FF_q$-linear spaces. Let $M_1=(E_1,\rho_1)$ and $M_2=(E_2,\rho_2)$ be a pair of lattice-equivalent $q$-polymatroids.
    Then $p(M_1;z) = p(M_2;z)$.
\end{Lemma}

\begin{proof}
    Let $\phi:\kL(E_1) \longrightarrow \kL(E_2)$ be a lattice isomorphism such that 
    $\rho_2(\phi(X)) = \rho_1(X)$ for all $X \in \kL(E_1)$.
    Since $\kL(E_1)$ and $\kL(E_2)$ are equivalent lattices, we have that $\dim(X) = \dim(\phi(X))$ for all $X \in \kL(E_1)$ and in particular $\mu_1(0,X) = \mu_2(0,\phi(X))$, where $\mu_i$ denotes the M\"{o}bius function on $\kL(E_i)$. Moreover, $X \leq Y$ in $\kL(E_2)$ if and only if $\phi(X) \leq \phi(Y)$ in $\kL(E_1)$. By assumption, $\ell_{M_1}(X) = \ell_{M_2}(\phi(X))$ for each $X \in \kL(E_1)$ and so the result now follows.
\end{proof}

We have the following results on the characteristic polynomial of the contraction of $M$ to a subspace $T$. These will be important later when we define the $q$-polymatroid version of the rank weight enumerator.

\begin{Lemma}\label{lem:ellM/Tperp}
	Let $T \leq E$ and $M=(E, \rho)$ a $q$-polymatroid. The following hold.
	\begin{enumerate}
		\item $\ell_{M.T}(X/T^\perp) = \ell_{M/T^\perp}(X/T^\perp) =\ell(X).$
		\item $\displaystyle p(M.T;z) = \sum_{T^\perp \leq X \leq E} \mu(T^\perp,X) z^{\ell(X)}$.
\item $\displaystyle p(M/T;z) = \sum_{T \leq X \leq E} \mu(T,X) z^{\ell(X)}$.	
	\end{enumerate}
\end{Lemma}
\begin{proof}
The first part follows from a direct computation:
\begin{align*}
\ell_{M.T}(X/T^\perp) =\ell_{M/T^\perp}(X/T^\perp) & = \rho_{M/T^\perp}(E/T^\perp)-\rho_{M/T^\perp}(X/T^\perp) \\
 & = \rho(E)-\rho(T^\perp) - \rho(X)+\rho(T^\perp) \\
 & = \rho(E)-\rho(X) = \ell(X).
\end{align*}
Let $\bar{\mu}$ denote the M\"obius function on the lattice of subspaces of $E/T$. Then, applying (1) we have:
\begin{align*}
	p(M.T;z) 
	&= p(M/T^\perp;z) = \sum_{T^\perp \leq X \leq E} \bar{\mu}(0,X/T^\perp) z^{\ell_{M/T^\perp}(X/T^\perp)}= \sum_{T^\perp \leq X \leq E} \mu(T^\perp,X) z^{\ell(X)},
\end{align*}
which proves (2). The last item follows directly from $M.T=M/T^\perp$ (Definition \ref{Cont}).
\end{proof}
Clearly, if $T$ has dimension $t$, then 
$\displaystyle p(M.T;z)	= \sum_{j=0}^t (-1)^j q^{\binom{j}{2}} \sum_{T^\perp \leq Y, \dim(Y)=n-t+j} z^{\ell(Y)}.$

\begin{Example}\label{exVamosContr}
We calculate $p(M.T;z)$ where $M$ is the V\'{a}mos $q$-matroid (Example \ref{exVamos}).
Let $T$ be a subspace of $E=\FF_q^8$.
If $T$ has  dimension $5$, then $\dim T^\perp=3$. We need only consider two cases, depending on whether or not $T^\perp$ is contained in a circuit (a member of $\mC$). Note that the circuits intersect pairwise in dimension $2$ or $0$, so $T^\perp$ cannot be in more than one circuit. \\
Suppose $T^\perp$ is in none of the circuits. Then for all $X$ such that $T^\perp< X\leq E$ we have that $\ell(X)=0$. For $X=T^\perp$, we have $\ell(X)=1$. So the $q$-matroid $M.T$ is lattice-equivalent to the uniform $q$-matroid $U_{1,5}$. Its characteristic polynomial is
\[ p(M.T;z)=\mu(T^\perp,T^\perp)z^1+\sum_{T^\perp< X\leq E}\mu(T^\perp,X)z^0 = z-1. \]
Suppose now that $T^\perp$ is contained in a circuit $C \in \mC$.
Among all $X$ such that $T^\perp\leq X\leq E$ we have that $\ell(X)=1$ for $X=T^\perp$ and $X=C$. Otherwise, $\ell(X)=0$. The $q$-matroid $M.T$ has rank $1$ and all $1$-dimensional spaces are independent, except for the circuit $C/T^\perp$. For the characteristic polynomial we get the following:
\begin{align*}
p(M.T;z) &=\mu(T^\perp,T^\perp)z+\mu(T^\perp,C)z + \sum_{{X\,:\,}T^\perp< X\leq E, X\neq C}\mu(T^\perp, X) \\
 &= z-z-\left(\qbin{8-3}{4-3}{q}-1\right)+q\qbin{8-3}{5-3}{q}-q^3\qbin{8-3}{6-3}{q}+q^6\qbin{8-3}{7-3}{q}-q^{10}, \\
 &= -\left(\qbin{5}{1}{q}-1\right)+q\qbin{5}{2}{q}-q^3\qbin{5}{3}{q}+q^6\qbin{5}{4}{q}-q^{10}=0.
\end{align*}
In fact, since $\sum_{{X\,:\,}T^\perp< X\leq E, X\neq C}\mu(T^\perp, X) +1 +\mu(T^\perp,C)=0$, we see that the constant term is zero without the need for any calculation.
\end{Example}

We continue to develop technical properties of the characteristic polynomial of the contraction $M.T$. In Section \ref{sec:designs}, we will use the fact that the characteristic polynomial of $M.T$ is identically zero when $T$ is an independent space of the dual $q$-polymatroid.

\begin{Lemma}\label{lem:indlperp0}
$T \leq E$ is an independent space of $M^*$ if and only if $\ell(T^\perp)=0$.
\end{Lemma}	
\begin{proof}
We have 
\begin{align*}
\ell(T^\perp) &= \rho(E)-\rho(T^\perp)
=\rho(E)-(\rho^*(T)+\rho(E)-r\dim T)
=r \dim(T)-\rho^*(T).
\end{align*}
Hence $T$ is an independent space of $M^*$ if and only if $\ell(T^{\perp})=0$.
\end{proof}

\begin{Lemma}\label{lem:pm.tneq0}
If  $T \leq E$, $\dim T >0$ is an independent space of $M^*${, then} $p(M.T;z)=0$.
\end{Lemma}
\begin{proof}
By Lemma \ref{lem:indlperp0}, $\ell(T^\perp)=0$. Since all subspaces of an independent space are independent, we have that $\ell(X)=0$ for all $T$ such that $T^\perp\leq X$. Applying this to the characteristic polynomial, we get
\[ p(M.T;z) = \sum_{T^\perp \leq X \leq E} \mu(T^\perp,X) z^{\ell(X)} = \sum_{T^\perp \leq X \leq E}\mu(T^\perp,X) =0. \qedhere \] 
\end{proof}

\begin{Lemma}\label{lem:pm.tcirc}
	Let $T \leq E$ be a circuit of $M^*=(E,\rho^*)$. 
	Then $p(M.T;z)=z^{\ell(T^\perp)}-1$.
\end{Lemma}
\begin{proof}
Let $X \leq E$. If $T^\perp$ is strictly contained in $X$, then $X^\perp$ is strictly contained in $T$,  and so $X^\perp$ is independent in $M^*$.  Therefore,   Lemma \ref{lem:indlperp0} gives that $\ell(X)=0$. Hence 
\begin{align*}
    p(M.T;z)&= \sum_{T^\perp \leq X \leq E} \mu(T^\perp,X) z^{\ell(X)} = z^{\ell(T^\perp)} + \sum_{T^\perp < X \leq E} \mu(T^\perp,X) \\
    &= z^{\ell(T^\perp)}-\mu(T^\perp,T^\perp) = z^{\ell(T^\perp)}-1. \qedhere
\end{align*} 
\end{proof}

\begin{Remark}\label{rem:cocirccpoly}
   Note that if $M$ is a $q$-matroid, a cocircuit $T$ of $M$ has $\ell(T^\perp)=\dim(T)-\rho^*(T)=\dim(T)-(\dim(T)-1)=1$ hence $p(M.T;z)=z-1$.
\end{Remark}

\begin{Lemma}\label{lem:ellresttperp}
	Let $M=(E,\rho)$ be a $q$-polymatroid and let $T \leq E$ be an independent space of $M^*$. The following hold.
	\begin{enumerate}
		\item $\rho(E)=\rho(T^\perp)$.
		\item For any subspace $U\leq T^\perp$, we have $\ell_{{M|T^\perp}}(U)=\ell(U).$
	\end{enumerate}
	
\end{Lemma}	

\begin{proof}
    By definition of the dual $q$-polymatroid, we have $\rho(T^\perp)=\rho^*(T)-{r}\dim(T)+\rho(E)$.
    Since $T$ is independent in $M^*$, $\rho^*(T)={r}\dim(T)$ and so we get $\rho(T^\perp)=\rho(E)$, which establishes (1).	
	Therefore, 
   $\ell_{{M|T^\perp}}(U) = \rho_{{M|T^\perp}}(T^\perp) - \rho_{{M|T^\perp}}(U)= \rho(T^\perp)-\rho(U)
   = \rho(E)-\rho(U)=\ell(U),$ which proves (2). 	
\end{proof}

\begin{Corollary}\label{cor:carpolyMTperp}
    Let $T\leq U$ be subspaces of $E$ such that $T$ is independent in $M^*$. If $U/T$ is a circuit in $M^*/T${, then}
    $$p((M^*/T)^*.(U/T);z)=p({M|T^\perp}.(U^\perp)^{\perp(T^\perp)};z)  = z^{\ell(U^\perp)}-1. $$       	
\end{Corollary}
\begin{proof}
	Recall {from Lemma \ref{lem:dualMTperp}} that $M^*/T \cong (M|T^\perp)^*$ (and hence $(M^*/T)^* \cong M|T^\perp$) under the map $\phi:A/T \mapsto (A^\perp)^{\perp(T^\perp)}$ for any $A \leq E$ with $T \leq A$. In particular, if $U/T$ is a circuit in $M^*/T${, then} $\phi(U/T)$ 
	is a circuit in $({M|T^\perp})^*$. Moreover $\phi(U/T)^{\perp(T^\perp)} = U^\perp$.
	From Lemmas \ref{lem:ellresttperp} and \ref{lem:pm.tcirc} we have
	\[ p((M^*/T)^*.(U/T);z) = p({M|T^\perp}.\phi(U/T);z) = z^{\ell_{{M|T^\perp}}(U^\perp)}-1=z^{\ell(U^\perp)}-1. \qedhere \] 
\end{proof}

The following result will be used in the proof of Corollary \ref{firstcorollary}.

\begin{Lemma}\label{lem:2.11}
	Let $W \leq E$ and let $T\leq W$ be an independent space of $M^*$. 
	Then 
	\[p({M|T^\perp}/W^\perp;z) = \sum_{{A\,:\,} A+T =W} p(M.A;z).\]	
\end{Lemma}

\begin{proof}
	By Lemmas \ref{lem:ellM/Tperp} and \ref{lem:ellresttperp}, we have 
	$\ell_{{M|T^\perp}/W^\perp} (U/W^\perp) = \ell(U)$ for any subspace $U$ satisfying $T \leq {U^\perp} \leq W$. 
	Since $\displaystyle p(M/U;z) = \sum_{{A\,:\,}U \leq A \leq E} \mu(U,A) z^{\ell(A)},$ by applying the M\"{o}bius inversion formula we have
	$\displaystyle z^{\ell(U)} = \sum_{{A\,:\,}U \leq A \leq E} p(M/A;z) $. 
	Therefore, we have
	\begin{eqnarray*}
		p({M|T^\perp}/W^\perp;z)&=& \sum_{{U\,:\,} W^\perp \leq U \leq T^\perp} \mu(W^\perp,U) z^{\ell(U)}\\
		 &=& \sum_{U{\,:\,} W^\perp \leq U \leq T^\perp} \mu(W^\perp,U) \sum_{{A\,:\,} U \leq A \leq E} p(M/A;z)\\
		 &=& \sum_{{U\,:\, }W^\perp \leq U \leq T^\perp} \mu(W^\perp,U) \sum_{{A\,:\, U} \leq A \leq E} p(M.A^\perp;z)\\
		 &=& \sum_{{V\,:\,} W^\perp \leq V^\perp \leq T^\perp} \mu(W^\perp,V^\perp) \sum_{A{\,:\,} V^\perp \leq A \leq E} p(M.A^\perp;z)\\	 
		 &=& \sum_{{V\,:\,} T \leq V \leq W} \mu(W^\perp,V^\perp) \sum_{A{\,:\,} 0 \leq A \leq V} p(M.A;z)\\
		 &=& \sum_{A{\,:\,} 0 \leq A \leq W} p(M.A;z) \sum_{V{\,:\,}A+ T \leq V \leq W} \mu(W^\perp,V^\perp) \\
		 &=& \sum_{A{\,:\,} A+T = W} p(M.A;z), 
	\end{eqnarray*}
    where the last equality follows from the fact that 
    \[\sum_{V{\,:\,}A+ T \leq V \leq W} \mu(W^\perp,V^\perp)=\sum_{{V\,:\,}W^\perp \leq V^\perp \leq A^\perp \cap T^\perp } \mu(W^\perp,V^\perp)
     = \left\{ \begin{array}{cl}
    		1 & \text{ if } A^\perp \cap T^\perp = W^\perp,\\
    		0 & \text{ otherwise. }
    	\end{array} \right.\qedhere
    \]  
   
\end{proof}

We now present some further results on the characteristic polynomial.

\begin{Definition}
A {\bf loop} of the $(q,r)$-polymatroid $M$ is a circuit of dimension $1$.
\end{Definition}

Clearly, if $e$ is a loop of $M$, then $0\leq \rho(e) < r$. If $M$ is a $q$-matroid{, then} the loops of $M$ all have rank zero.

\begin{Lemma}\label{lem:isthmus}
	Let $e$ be a one-dimensional subspace of $E$. The following are equivalent:
	\begin{enumerate}
		\item $p(M.e;z)=0$,
		\item $\rho(e^\perp) = \rho(E)$,
		\item $e$ is not a loop in $M^*$.
	\end{enumerate}
\end{Lemma}

\begin{proof}
	We have $p(M.e;z)=z^{\ell(e^\perp)}-z^{\ell(E)}=z^{\ell(e^\perp)}-1$, which is zero if and only if 
	$\ell(e^\perp) = \rho(E)-\rho(e^\perp) = 0$. This shows that (1) and (2) are equivalent.
	The one-dimensional space $e$ is independent in $M^*$ if and only if $\rho^*(e)=r$.
	Since $\rho^*(e)=r\dim(e)-\rho(E)+\rho(e^\perp)=r-\rho(E)+\rho(e^\perp)$, this occurs if and only if 
	$\rho(e^\perp)=\rho(E),$ which shows that (2) and (3) are equivalent.
\end{proof}

\begin{Remark}
We remark that for a classical matroid $M$, if $e$ is a loop of $M${, then} it is easy to see that $p(M.e;z)=0$. This is because if $e$ is a loop, by semimodularity, the fact that $e\notin E-e$ forces $\rho(E-e)=\rho(E)$. Then by the classical version of the lemma above, this gives that $p(M.e;z)=0$. For a $q$-matroid, however, it may occur that $e \leq e^\perp$, in which case if $e$ is a loop, the same argument using semimodularity does not imply that $\rho(e^\perp)=\rho(E)$.
\end{Remark}

\begin{Definition}
   For each $A \in \kL(E)$, define $c(A):=\{ X \leq E : A \leq X,\rho(A) =\rho(X) \}$.
   The {\bf closure} of $A$ in the $(q,r)$-polymatroid $M$ is denoted by $\cl(A)$ and is defined to be the vector space sum of the members of $c(A)$; that is, $\cl(A) := \sum_{ X \in c(A)} X$.
\end{Definition}

\begin{Lemma}\label{lem:monic}
	Let $L=\cl(\{0\})$. Let $X$ be a subspace of $E$ such that $X^\perp \leq L$. Then
	\begin{eqnarray*}
p(M.X;z)
	&=&\left\{ 
    	\begin{array}{rl}
    	    \displaystyle z^{\rho(E)}+\sum_{{A\,:\,}X^\perp\lneq A \leq E} \mu(X^\perp,A) z^{\ell(A)}  
    		& \text{ if }X = L^\perp ,\\
    		\displaystyle \sum_{{A\,:\,}X^\perp\leq A \leq E, A \nleq L} \mu(X^\perp,A) z^{\ell(A)}
    		& \text{ otherwise. } 
    	\end{array}
    	\right.   
    \end{eqnarray*}
    If $X^\perp = L${, then} $p(M.X;z)$ is a monic polynomial of degree $\rho(E)$ in $z$.
	In particular, if $M$ is a $q$-matroid and has no loops{, then} $p(M;z)$ is monic polynomial of degree $\rho(E)$. 
\end{Lemma}

\begin{proof}
    From Lemma \ref{lem:ellM/Tperp} we have that
    \begin{eqnarray*}
  p(M.X;z)=    \sum_{{A\,:\,}X^\perp\leq A \leq E} \mu(X^\perp,A) z^{\ell(A)}&=& \\
= z^{\rho(E)} \sum_{{A\,:\,}X^\perp\leq A \leq L} \mu(X^\perp,A) &+& \sum_{A\,:\,X^\perp\leq A \leq E, A \nleq L} \mu(X^\perp,A) z^{\ell(A)}. \\
    \end{eqnarray*}
    By the definition of the M\"obius function, $\sum_{{A\,:\,X^\perp}\leq A \leq L} \mu(X^\perp,A)=0$
    unless $X^\perp = L$.
    If $A \nleq L${, then} $\ell(A)=\rho(E)-\rho(A)<\rho(E)$, so if $L=X^\perp$, then $p(M.X;z)$ is a monic polynomial with leading term $z^{\rho(E)}$.
     Furthermore, setting $X=E$, we obtain that if $M$ is a $q$-matroid with no loops, then its characteristic polynomial is monic of degree $\rho(E)$.
\end{proof}

\begin{Remark}
If $M$ is a $q$-matroid, $\cl(\{0\})$ is the space containing all the loops. 
\end{Remark}

In the $q$-matroid case, cryptomorphisms  between  axiom systems such as those relating to independent spaces, the closure function, flats, hyperplanes etc., were established in \cite{bcj}. {We therefore have the following facts, as in the case for classical matroids. The reader is referred to \cite{bcj} and the references therein for further details.
A subspace $F$ is called a {\bf flat} of a $q$-matroid if $\cl(F)=F$. For each $B \leq E$, there is a unique flat $F$ such that $\cl(B)=F$, in which case $\rho(B) = \rho(F)$. Moreover, if $M$ is a $q$-matroid, its collection of flats forms a semi-modular lattice \cite{WINEpaper1}.
A {\bf hyperplane} $H < E$ is a flat that is maximal with respect to containment, that is,
if $H \leq F$ for some flat $F${, then} either $F=E$ or $F=H$. 
Every flat of $M$ is an intersection of hyperplanes and  for every  hyperplane $H$, we have that $H^\perp$ is a cocircuit of $M$. Therefore,  for every flat $F$ of $M$, $F^\perp$ is the vector space sum of a collection of cocircuits.
}

The following result will be used in Corollary \ref{cor:pcocirc}.

\begin{Theorem}\label{th:cocircs}
    Let $M$ be a $q$-matroid. Let $X \leq E$ contain a unique cocircuit $C$.  
    Then
	\begin{eqnarray*}
p(M.X;z)
	&=&\left\{ 
    	\begin{array}{ll}
    	    \displaystyle z-1
    		& \text{ if } X = C,\\
    		0
    		& \text{ otherwise. } 
    	\end{array}
    	\right.   
    \end{eqnarray*}
\end{Theorem}

\begin{proof}
   If $X=C$, then by Remark \ref{rem:cocirccpoly} we have that $p(M.X;z)=z-1$. Assume now that $C \lneq X$. Then $X$ is not a sum of cocircuits of $M$ and hence $X^\perp$ is not a flat. 
    Clearly $X$ is a dependent space of $M^*$
    and by the uniqueness of $C$, any subspace of $X$ that is dependent in $M^*$ contains $C$. 
    Therefore, by Lemma \ref{lem:indlperp0}, $\ell(A^\perp)=0$ for every $A \leq X$ such that $C \nleq A$.
    
   Let $F$ be a flat of $M$.
    For any $A \leq E$ such that $\cl(A^\perp)=F$,  we have  
    $\rho(A^\perp)=\rho(F)$ and hence $\ell(A^\perp)=\ell(F)$.
  Furthermore, $\text{cl}(C^\perp) = C^\perp$ since $C$ is a cocircuit of $M$.
    This will be used in the following computation
    of $p(M.X;z)$: 
        \begin{align*}
        p(M.X;z) & = \sum_{{A\,:\,} X^\perp \leq A \leq E} \mu(X^\perp,A) z^{\ell(A)} 
                   =  \sum_{{A\,:\,} A \leq X } \mu(A,X) z^{\ell(A^\perp)}\\
                 & = \sum_{{A\,:\,} C  \leq A \leq X } \mu(A,X) z^{\ell(A^\perp)} +
                     \sum_{A{\,:\, A }\leq X, C \nleq A } \mu(A,X)\\ 
                & = \sum_{\substack{{F\,:\,}\cl(F) =F \\ X^\perp \leq F \leq C^\perp }}\sum_{\substack{{A\,:\, X^\perp} \leq A^\perp \leq C^\perp \\ \cl(A^\perp)=F }} \mu(X^\perp,A^\perp) z^{\ell(F)} +
                     {\sum_{A: A \leq X} \mu(A,X) -\sum_{A: C \leq A \leq X} \mu(A,X)}\\     
                 & = \sum_{\substack{{F\,:\,}\cl(F) =F \\ X^\perp \leq F \leq C^\perp }}z^{\ell(F)}\sum_{\substack{{A\,:\,} X^\perp \leq A^\perp \leq C^\perp \\ \cl(A^\perp)=F }} \mu(X^\perp,A^\perp)
                     .\\    
    \end{align*}  
    Since the lattice of flats of $M$ is semi-modular, 
    by \cite[Proposition 3.3]{whittle}, we have
    \[ \sum_{\substack{{A\,:\,} X^\perp \leq A^\perp \leq C^\perp\\ \cl(A^\perp)=F} } \mu(X^\perp,A^\perp) = 0, \]
    and so the result follows.
\end{proof}

\begin{Remark}
   In fact, by a similar argument (also essentially the same as for classical matroids),  for a $q$-matroid $M$,  we have  $p(M.X;z) = 0$ unless $X^\perp$ is a flat in $M$. 
   Equivalently, we have that $p(M.X;z)=0$ unless $X$ is a sum of cocircuits of $M$.
\end{Remark}

\subsection{The Weight Enumerator of a {\em q}-Polymatroid}

We  next define the {\em weight enumerator} of a $q$-polymatroid. In Section \ref{sec:macwill}, we will show that  its values satisfy a duality property and in Section \ref{sec:designs}, we will apply this duality result to establish a criterion for identifying a weighted subspace design determined by a $q$-polymatroid.

\begin{Definition}\label{def:weighten}
	We define the {\bf weight enumerator} of the $(q,r)$-polymatroid $M$ to be the list $[A_M(i;z) : 0 \leq i \leq n]$, where for each $i$ we define
	$$A_M(i;z):=\sum_{X \leq E, \dim(X)=i} p(M.X;z)= 
	\sum_{X \leq E, \dim(X)=i} p(M/X^\perp;z) .$$
\end{Definition}

\begin{Lemma}\label{lem:M/T}
	Let $T\leq E$. The following hold.
	\begin{enumerate}
		\item If  $T \leq Z \leq E${, then} $p((M/T)\big/(Z/T);z) = p(M.Z^\perp;z)$.
		\item $\displaystyle A_{M/T}(j;z) = \sum_{X \leq T^\perp: \dim(X)=j} p(M.X;z)$.
	\end{enumerate}	
\end{Lemma}
\begin{proof}
Let $T\leq Z\leq Y\leq E$. Then $(Y/T)\big/(Z/T)$ and $Y/Z$ are isomorphic.
Let $V=(E/T)\big/ (Z/T)$ and write $M_V=(M/T)\big/ (Z/T)$. We have a lattice isomorphism between $\mL(E/Z)$ and $\mL(V)$. Moreover, it is easy to check that $\rho_{M_V}((Y/T)\big/(Z/T)) = \rho_{M/Z}(Y/Z)$. Therefore, $M_V$ and $M/Z$ are lattice-equivalent $q$-polymatroids.
We thus have
	\[	
	    p(M_V;z) =  p(M/Z;z) = p(M.Z^\perp;z) . 
    \] 
    Let $X \leq T^\perp$. It is straightforward to check that $\dim((X^\perp/T)^{\perp(E/T)})=\dim(X)$. Therefore,
    \begin{align*}
    	A_{M/T}(j;z)&= \sum\limits_{\substack{X: X^\perp/T\leq E/T,\\ \dim((X^\perp/T)^{\perp(E/T)})=j}} p((M/T) \big/ (X^\perp/T);z) = \sum \limits_{\substack{X: X \leq T^\perp \\ \dim(X)=j}} p(M.X;z). \qedhere
    \end{align*}
\end{proof}

\section{Matrix Codes and {\em q}-Polymatroids}\label{sec:codes}
	
	We consider properties of a $q$-polymatroid arising from an $\FF_q$-linear rank-metric code.
	There are several papers outlining properties of rank-metric codes.
	The $q$-polymatroids associated with these structures have been studied in 
	\cite{GLJ,gorla2019rank,shiromoto19}. 

   	\begin{Notation}
   	Throughout this section, we let $m$ be a positive integer {and $E=\mathbb{F}_q^n$}. 
   	We write $U^\perp$ to denote the orthogonal complement of $U \leq E$ with respect to a non-degenerate symmetric bilinear form $b_E$ on $E$. By abuse of notation, we also write $U^\perp$ to denote the orthogonal complement of 
   	\begin{itemize}
   	    \item $U \leq \FF_q^{n \times m}$ with respect to the inner product
	$b_{\FF_q^{n \times m}}$ defined by $b_{\FF_q^{n \times m}}(X,Y) = Tr(XY^T)$ for all $X,Y \in \FF_q^{n \times m}$ and
	\item $U \leq \FF_{q^m}^n$ with respect to the dot product defined by 
	$x \cdot y = \sum_{i=1}^n x_i y_i$ for all 
	$x=(x_1,\ldots,x_n), y=(y_1,\ldots y_n) \in \FF_{q^m}^n$.
   	\end{itemize}
   	\end{Notation}

	\begin{Definition}
	    We say that $C \subseteq  \FF_q^{n \times m}$ is a
	    {\bf linear rank-metric code}, or a {\bf matrix code} if $C$ is a subspace
	    of $\FF_{q}^{n \times m}$. The {\bf minimum distance} of $C$ is the minimum rank of any 
	    {nonzero} member of 
	    $C$. We say that $C$ is an $\FF_q$-$[n \times m, k,d]$ rank-metric code if it has $\FF_q$-dimension $k$ and minimum distance $d$. 
	    The {\bf dual code} of $C$ is 
	    $C^\perp :=\{ Y \in \FF_q^{n \times m} : {\rm Tr}(XY^T) =0 \: \forall\: X \in C  \}$.
	    {Finally, for each $i \in \{0,\ldots,n\}$, we define $W_i(C) := |\{ A \in C : \rk(A) = i \}|$. The list $[W_i(C) : 0 \leq i \leq n]$ is called the {\bf weight distribution} of $C$.}
	\end{Definition}

	\noindent{For $X \leq E$ we write $\col(X)$ to denote the column space of $X$ over $\FF_q$.}
	
	\begin{Definition}
	Let $X \in \FF_q^{n \times m}$ and let $U\leq E$. 
	We say that $U$ is the support of $X$ if $\col(X)=U$.
	Let $C$ be an $\FF_q$-$[n \times m, k,d]$ rank-metric code.
	We say that $U$ is a {\bf support} of $C$ if there exists some $X \in C$
	with support $U$.
	\end{Definition}

	\begin{Definition}\label{def:codepoly}
	   Let $m$ be a positive integer and let $C$ be an $\mathbb{F}_{q}$-$[n\times m,k,d]$ rank-metric code.
	   For each subspace 
	   $U \leq E$, we define
	   $$C_U:=\{A \in C : \col(A) \leq U^\perp\} \text{ and }
	   C_{=U}:=\{A \in C: \col(A) = U^\perp\}.$$ 
	   Let $\rho: \mL (E) \longrightarrow \NN_{\geq 0}$ be defined by
	   	$\displaystyle \rho(U):=k-\dim(C_U).$
	   	{Then $(E,\rho)$ is a $(q,m)$-polymatroid \cite[Theorem 5.3]{gorla2019rank} and we denote it by $M_C$.}	
	\end{Definition}
	Clearly, we have $\displaystyle \ell(U) = \dim(C_U)$ 
	for every $U \leq E$.  
	
	\begin{Lemma}\label{lem:repqpoly}
		Let $C$ be an $\mathbb{F}_{q}$-$[n\times m,k,d]$ rank-metric code.
		The following hold.
		\begin{enumerate}
			\item $M_{C^\perp} = (M_C)^*$.
  			\item $p(M_C/U;q)=|C_{=U}|$. 
			\item $W_i(C) = A_{M_C}(i;q)$ for each $i \in \{1,\ldots,n\}$.
			\item  $A_{M_C}(i;q)=0$ if and only if $p(M_{C}/U;q) =0$ for every $i$-dimensional subspace $U \leq E$.
			\item If  $A_{M_C}(i;q)=0$, then $A_{M_C/T}(i;q)=0$ for every 
			subspace $T \leq E$.  
		\end{enumerate} 
	\end{Lemma}
	
	\begin{proof}
		(1) has been established in \cite[Theorem 7.1]{gorla2019rank}.	
		Let $M=M_C$.
		Since $\displaystyle |C_U| = \sum_{V: U \leq V} |C_{=V}|,$ by M\"obius inversion we have
		$$|C_{=U}| = \sum_{V: U \leq V} \mu(U , V)|C_V| = \sum_{V: U \leq V} \mu(U , V)q^{\ell(V)} = p(M.U^\perp;q)=p(M/U;q).$$
		Therefore (2) holds. 
		The number of codewords of $C$ that have rank $i$ over $\mathbb{F}_q$ is
		$$W_i(C)=\sum_{{U\,:\,}\dim(U)=n-i} |C_{=U}| =\sum_{{U\,:\,}\dim(U)=n-i} p(M.U^\perp;q) = \sum_{{U\,:\,}\dim(U)=i} p(M.U;q)=A_{M}(i;q),$$
		which gives (3). 
		Clearly, $A_{M}(i;q) = 0$ if and only if $p(M.U;q) = 0$ for each $U\leq E$ of dimension $i$, which gives (4).
		Let $T$ be a 
		subspace of $E$. By Lemma \ref{lem:M/T} we have
		$$A_{M/T}(i;q)= \sum_{X \leq {T^\perp\,:\, \dim(X)=i}} p(M.X;q) .$$
		If $A_{M}(i;q) = 0$, then from (4) we have $p(M.X;q)=0$ for each $i$-dimensional subspace $X$, and so we get  $A_{M/T}(i;q)=0$, which proves (5).
	\end{proof}
	
 \begin{Remark}
    Note that Part (2) of Lemma \ref{lem:repqpoly} is an instance of the {Critical Theorem~\cite{CrapoRota}} for $q$-polymatroids and matrix codes. 
 \end{Remark}
 
  \begin{Remark}
     In \cite{gorla2019rank}, the authors define a pair of $q$-polymatroids associated
     with a matrix code. The one given above is the $q$-polymatroid whose rank function is
     determined by the column-spaces of the codewords. A second $q$-polymatroid is  the one whose rank function is determined by the row spaces of the codewords.
  \end{Remark}
  
  One way to construct an $\FF_q$-$[n \times m,k,d]$ rank-metric code is by taking a subspace of $\FF_{q^m}^n$,
  and expanding its elements with respect to a basis of $\FF_{q^m}$ over $\FF_q$. Such rank-metric codes are referred to as vector rank-metric codes.
  
  \begin{Definition}
	Let $\Gamma$ be a basis of $\FF_{q^m}$ over $\FF_q$. 
	For each $x \in \FF_{q^m}^n$, we write $\Gamma(x)$ to denote the $n \times m$ matrix over $\FF_q$ whose $i$th row is the coordinate vector of
	the $i$th coefficient of $x$ with respect to the basis $\Gamma$.
	The \textbf{rank} of $x$ is the rank of the matrix $\Gamma(x)$. Note that the rank of $x$ is well-defined, being independent of the choice of basis $\Gamma$.
\end{Definition}

For the remainder, we fix $\Gamma$ to be a basis of $\FF_{q^m}$ over $\FF_q$.

\begin{Definition}
	A (linear rank-metric) \textbf{vector code} $C$ is an $\FF_{q^m}$-subspace of $\FF_{q^m}^n$.
	The \textbf{minimum distance} of $C$ is the minimum rank of any non-zero element of $C$.
	We say that $C$ is an $\FF_{q^m}$-$[n,k,d]$ code if it has ${\FF_{q^m}}$-dimension $k$ and $\Gamma(C)$ has minimum rank distance $d$.
	The code $C^\perp$ denotes the \textbf{dual code} of $C$ with respect to the standard dot product on $\FF_{q^m}^n$. 
\end{Definition}
 
Each vector rank-metric code determines a $q$-matroid, as follows.
	
	\begin{Definition}
	   Let $C$ be an $\mathbb{F}_{q^m}$-$[n,k,d]$ rank-metric code.
	   Let $U \leq E$ and let $x \in C$. We say that $U$ is a {\bf support} of $x$ if
	   $U$ is the column space of $\Gamma(x)$ and we write $\sigma(x) = U$.
	   For each subspace 
	   $U \leq E$, 
	   we define
	   $$C_U:=\{x \in C : \sigma(x) \leq U^\perp\} \text{ and }
	   C_{=U}:=\{x \in C: \sigma(x) = U^\perp\}.$$ 
	   Let $\rho: \mL (E) \longrightarrow \NN_{\geq 0}$ be defined by
	   	$\displaystyle \rho(U):=k-\dim_{\Fqm}(C_U).$
	   {Then $(E,\rho)$ is a $q$-matroid \cite[Theorem 24]{JP18} and we denote it by $M_C$.} 	
	\end{Definition}

	\begin{Remark}\label{rem:mols}
	Note that in the definition given above, the rank function for the $q$-matroid
	of $C$ is the rank function of the associated $(q,m)$-polymatroid as defined in Definition \ref{def:codepoly}, divided by $m$. Since $C$ is $\FF_{q^m}$-linear, $C_U$ is an $\FF_{q^m}$-vector space for each subspace $U$ and so has $\FF_q$-dimension a multiple of $m$.
	Therefore the results of Lemma \ref{lem:repqpoly} hold with $q^m$ in place of $q$. 
	For example, with respect to the characteristic polynomial of the $q$-matroid, we have $p(M/U;q^m)=|C_{=U}|$ for an $\FF_{q^m}$-$[n,k,d]$ code $C$ and subspace $U$. 
	\end{Remark}
	
	Let $C$ be an $\FF_{q^m}$-$[n,k,d]$ code. 
	Recall that for any $U \leq \FF_{q}^n$ we have 
        \[
          \ell_{ M_C}(U) = \dim_{\FF_{q^m}}(C_U) = \dim(\{ x \in C : \sigma(x) \leq U^\perp\}).
        \]
        Now $U^\perp$ is independent in $M_{C^\perp}$ if and only if $\ell_{ M_C}(U)=0$, 
        which occurs if and only if no subspace of $U^\perp$ is a support of $C$.
        Therefore every support of $C$ corresponds to a dependent space of $M_{C^\perp}$. 
        
 In the following example we illustrate the notions discussed in Sections \ref{sec:charpoly} and \ref{sec:codes}. We calculate the characteristic polynomial of $M_C$ by carefully studying the structure of the $q$-matroid and its dual.
	
	\begin{Example}\label{ex:code}
	    {Let $\alpha$ be a root of $x^6+x^4+x^3+x+1 \in \FF_{2^6}[x]$. Then $\alpha$ is a primitive element of $\FF_{2^6}$.} Let $C$ be the $\FF_{2^6}$-$[6,3,3]$ vector rank-metric code generated by the matrix:
	    \[
	    G = 
	    \left[
	    \begin{array}{cccccc}
	        1 & 0 & 0 & \alpha^{13} & \alpha^{47} & \alpha^{35} \\
	        0 & 1 & 0 & \alpha^{44} & \alpha^{62} & \alpha^{32} \\
	        0 & 0 & 1 & \alpha^{34} & \alpha^{22} & \alpha^{19}
	    \end{array}
	    \right].
	    \]
	    With respect to the basis $\Gamma=\{1,\ldots,\alpha^5\}$, the rows of $G$ are expanded to the following binary matrices: 
	    \[
	     \left[
	    \begin{array}{cccccc}
	        1 & 0 & 0 & 0 & 0& 0\\
            0 & 0 & 0 & 0 & 0& 0\\
            0 & 0 & 0 & 0 & 0& 0\\
            1 & 0 & 0 & 1 & 0& 0\\
            0 & 1 & 0 & 1 & 1& 1\\
            1 & 0 & 0 & 0 & 1& 0\\
	    \end{array}
	    \right],
	     \left[
	    \begin{array}{cccccc}
	        0& 0& 0& 0& 0& 0\\
            1& 0& 0& 0& 0& 0\\
            0& 0& 0& 0& 0& 0\\
            0& 0& 1& 1& 1& 1\\
            1& 0& 1& 1& 0& 1\\
            0& 1& 0& 0& 1& 1\\
	    \end{array}
	    \right],
	     \left[
	    \begin{array}{cccccc}
	       0& 0& 0& 0& 0& 0\\
           0& 0& 0& 0& 0& 0\\
           1& 0& 0& 0& 0& 0\\
           1& 0& 1& 0& 0& 1\\
           0& 0& 1& 1& 1& 0\\
           0& 1& 1& 1& 0& 1\\ 
	    \end{array}
	    \right].
	    \]
	    A basis of $\Gamma(C)$ over $\FF_2$, which has 18 elements, is found by multiplying each row of $G$ by successive powers of $\alpha$ and expanding with respect to $\Gamma$.  We have that 
	    $\Gamma(C)$ is an $\FF_2$-$[6 \times 6,18,3]$ rank-metric code with rank-metric weight distribution $[1,0,0,567,37044,142884,81648]$. Moreover, $C$ is formally self-dual, that is, its dual code has the same weight distribution as $C$.
	    Now consider the $q$-matroid  $M:=M_C$ arising from $C$, with rank function satisfying
	    $\rho(U)=3 - \dim_{\FF_{2^6}}(C_U)$ for each $U \leq \FF_2^6$.
        In Figure \ref{fig:ex} we write down the number of subspaces of $\FF^6_{2}$ for each possible value of $\ell(U)= \dim_{\FF_{2^6}}(C_U)$. 
        
        \begin{center}
     \begin{figure}[h!]
     \resizebox{\textwidth}{!}{
     \begin{tabular}{|c|c|c|c|c|c|c|c|}
     \hline
        \diagbox{{$\ell(U)$}}{$\dim(U)$} &$0$ & $1$ & $2$ & $3$ & $4$ & $5$ & $6$  \\
        \hline
        $0$ & ${0}$ & $0$ & $0$   & $\qbin{6}{3}{2}-9 = 1386 $ & $\qbin{6}{4}{2}=651$ & $\qbin{6}{5}{2}=63$ & $1$  \\
        \hline
        $1$ & $0$ & $0$ &  $\qbin{6}{2}{2}=651$ & $9$    & $0$              & $0$ & $0$ \\
        \hline
        $2$ & $0$ & $\qbin{6}{1}{2}=63$  & $0$ & $0$  & $0$ & $0$ & $0$\\
        \hline
        $3$ & $1$ & $0$  & $0$ & $0$  & $0$ & $0$ & ${0}$\\
     \hline
     \end{tabular} }
     \caption{{Number of subspaces for each possible $\ell(U)$ value.}}
     \label{fig:ex}
     \end{figure}
     \end{center}
     Using the entries of the table of Figure \ref{fig:ex}, we write down the characteristic polynomial of $M$:
     \begin{align*}
         p(M;z) & = \sum_{{U\, :\,} 0 \leq U \leq E} \mu(0,U) z^{\ell(U)}\\
                  & = z^3 + \sum_{{U \, :\, 0 \lneq} U \leq E, \ell(U)=2} \mu(0,U) z^2
                      +\sum_{{U \,:\, 0\lneq} U \leq E, \ell(U)=1} \mu(0,U) z 
                      + \sum_{U {\, : 0 \lneq U \leq E}, \ell(U)=0} \mu(0,U)\\
                & = z^3 -63z^2+1230z-1168
                =(z-1)(z^2-62z+1168).  
     \end{align*}
     
     We will explain the values in this table column by column from right to left: recall that to say something about $\ell(U)$, we have to consider how the supports of $C$ relate to $U^\perp$. Along the way, we will compute the different possible values
        of $p(M/U;q^m)=|C_{=U}|$, which, by Lemma \ref{lem:repqpoly}, counts the number of codewords of $C$ with support equal to $U^\perp$.
            
    Since the rank distance of $C$ is $3$, $C$ has no supports of dimension less than $3$ apart from $\{0\}$. Hence $\ell(U) = \dim_{\FF^6_{2}}(C_U) = 0$ for each of the 651 subspaces $U \leq \FF_{2}^6$ of dimension $4$ and each of the $63$ spaces of dimension $5$. So by Lemma \ref{lem:indlperp0}, the respective $1$ and $2$-dimensional orthogonal complements of these spaces are independent in $M^*$. We remark that by Lemma \ref{lem:repqpoly}, we have $M^*=M_{C^\perp}$. 
   
        We now consider the $3$-dimensional subspaces. 
        Since any proper subspace of a $3$-dimensional subspace $U\leq\FF_2^6$
        is independent in $M^*$, it must be the case that if $U$ is dependent in $M^*$, it is a cocircuit of $M$.
        Then by Remark \ref{rem:cocirccpoly}, we have $\ell(U^\perp)=1$ (indeed $p(M.U;z) = z-1$) and
        $p(M.U;2^6) = |C_{=U^\perp}| = 2^6-1 = 63$.
        Therefore, by inspection of the weight enumerator, we see that there are $9=567/63$, different $3$-dimensional spaces that are supports of $C$. We list the $3$-dimensional cocircuits of $M$ below:
        \begin{align*}
           \langle(0 1 0 0 1 1), (0 0 1 0 1 0), (0 0 0 1 0 0) \rangle,
           \langle (1 0 1 1 0 0), (0 1 0 0 0 0),(0 0 0 0 0 1)\rangle,
           \langle  (1 0 0 0 0 1), (0 1 1 0 0 0), (0 0 0 0 1 0)\rangle,\\
           \langle (1 0 0 1 1 1), (0 1 0 0 1 0), (0 0 1 1 0 1)\rangle,
           \langle (1 0 0 1 1 0), (0 1 0 1 0 1), (0 0 1 0 0 1)\rangle,
           \langle (1 0 0 0 1 0), (0 0 1 0 1 1), (0 0 0 1 1 1)\rangle,\\
           \langle (1 1 0 0 0 1), (0 0 0 1 0 1), (0 0 0 0 1 1)\rangle,
           \langle (1 0 0 1 0 0), (0 1 0 1 0 0), (0 0 1 1 1 1)\rangle,
           \langle (1 0 0 0 0 0), (0 1 0 1 1 0), (0 0 1 0 0 0)\rangle.
        \end{align*}
        
        Every other $3$-dimensional subspace $U$ is a non-support of $C$, as are all its nontrivial subspaces, hence $\ell(U^\perp)=\dim_{\FF_{2^6}}( C_{U^\perp})=0$. We remark that Lemma \ref{lem:indlperp0} says for such $U$ that $U^\perp$ is independent in $M^*$, and Lemma \ref{lem:pm.tneq0} gives that $p(M.U;z)=0$. 
         By computation, we obtain that there are $588$ $4$-dimensional supports of $C$ and that none of these spaces contains a cocircuit of dimension $3$. 
     Therefore, each such subspace $U$ is itself a cocircuit and so we have that $\ell(U^\perp)=1$, $p(M.U;z)=z-1$, and $p(M.U;2^6)=|C_{={U^\perp}}| =2^6{-1}$.
     There remain a further $\qbin{6}{4}{2}-588=63$ $4$-dimensional subspaces that are not supports of $C$.
     Every $3$-dimensional cocircuit is contained in $\qbin{6-3}{4-3}{2} =\qbin{3}{1}{2}= 7$ different $4$-dimensional spaces and every pair of $3$-dimensional cocircuits span $\FF_2^6$.
     Therefore, every $4$-dimensional non-support of $C$ contains at most one 
     $3$-dimensional cocircuit and since there are $9\cdot 7 = 63$ such $4$-dimensional non-supports altogether, 
     each of them contains a unique $3$-dimensional cocircuit. 
     It follows that $\ell(U^\perp)=1$ for every $4$-dimensional subspace $U$.
     
     By direct computation it can be checked that there are $63$ $5$-dimensional supports of $C$ and of course the only $6$-dimensional support is the entire space $\FF_2^6$.
     Each $5$-dimensional support $U$ is the support of exactly $2268$ different codewords, so 
     $p(M.U;2^6)=|C_{={U^\perp}}| = 2268$.
     Therefore, $\ell(U^\perp)=\dim_{\FF_{2^6}}(C_{U^\perp})\geq 2$. If $\ell(U^\perp)=3$ then the support of every codeword is contained in $U$, which is impossible as $C$ has words of rank 6. It follows that $\ell(U^\perp)=2$.  
     All computations carried out in this example were done using MAGMA~\cite{Magma}. 
	\end{Example}

\section{MacWilliams Identities for {\em q}-Polymatroids}\label{sec:macwill}

We establish a version of the MacWilliams identities for the $(q,r)$-polymatroids that we shall use in establishing criteria for the existence of a weighted $t$-design over $\Fq$. 
Duality via the rank polynomial of a $q$-polymatroid was considered in \cite{shiromoto19}.
We start with a result that relates the characteristic polynomial of a $q$-polymatroid to that of its dual. The statements of Theorem \ref{lem:charMM^*} and Corollary \ref{cor:charpolydual} may be regarded as $q$-analogues of \cite[Corollary 12]{britzshiromotomacwill}. However, unlike the proofs given here, which rely only on M\"obius inversion, the proof of \cite[Corollary 12]{britzshiromotomacwill} relies on an existing version of MacWilliams duality theorem for matroids, which shows that the weight enumerator polynomial of the dual of a matroid can be retrieved from the weight enumerator polynomial of the original matroid by a substitution of variables. A $q$-analogue of this result is not known.

\begin{Theorem}\label{lem:charMM^*}
	Let $U \leq E$. Then
	$$\sum_{A:A \leq U} p(M^*.A;z) = z^{r\dim(U)-{\rho(E)}} \sum_{A:A \leq U^\perp} p(M.A;z).$$
\end{Theorem}

\begin{proof}
We have by Lemma \ref{lem:ellM/Tperp} and then replacing $X$ with $A^\perp$ that
\[ p(M^*.U;z) = \sum_{X\,:\,U^\perp \leq X \leq E} \mu(U^\perp,X) z^{\ell^*(X)} =\sum_{A\,:\,A\leq U} \mu(A,U) z^{\ell^*(A^\perp)}. \]
To this we apply M\"{o}bius inversion (Lemma \ref{inversionformula} part (2)), duality (Definitions \ref{dualit} and \ref{funzionel}) and M\"{o}bius inversion again to get
\begin{align*}
\sum_{A\,:\,U\leq A \leq E} p(M^*.A;z) & = z^{\ell^*(U)} \\
& = z^{\ell(U^{\perp})-\rho(E)+r\dim(U^\perp)}\\
& = z^{r\dim(U^\perp)-\rho(E)} \sum_{A\,:\,U^\perp \leq A \leq E} p(M.A;z). \qedhere
\end{align*}

\end{proof}

We now show that for any subspace $U \leq E$, the characteristic polynomial of $M^*.U$ is completely determined by the set 
$\{(p(M.V;z), V): V \leq E \}$.

\begin{Corollary}\label{cor:charpolydual}
	Let $U \leq E$. We have the identity:
	\begin{equation*}
	z^{\rho(E)} p(M^*.U;z)  = \sum_{V \leq E} p(M.V;z) \sum_{j=0}^{\dim(U\cap V^\perp)} \qbin{\dim(U\cap V^\perp)}{j}{q}(-1)^{\dim(U)-j} q^{\binom{\dim(U)-j}{2}}z^{jr}.
	\end{equation*}
\end{Corollary}
\begin{proof}	
From Lemma \ref{lem:charMM^*}, we have:
\begin{equation}\label{eq:charMM^*}
  \sum_{{A\,:\,A} \leq U} p(M^*.A;z) = z^{r\dim(U)-{\rho(E)}} \sum_{{V\,:\,V} \leq U^\perp} p(M.V;z) .
\end{equation}
Apply the M\"{o}bius inversion formula to Equation~\eqref{eq:charMM^*}  to get the identity
$$\displaystyle p(M^*.U;z) = \sum_{{A\,:\, A} \leq U} \mu(A,U) z^{r\dim(A)-{\rho(E)}} \sum_{{V\,: V }\leq A^\perp} p(M.V;z).$$
Then
\begin{align*}
  z^{\rho(E)} p(M^*.U;z) & = \sum_{{A\,:\, A} \leq U} \mu(A,U) z^{r\dim(A)} \sum_{V \leq A^\perp} p(M.V;z)\\
   & = \sum_{V \leq E}p(M.V;z) \sum_{{A\,:\, A} \leq U\cap V^\perp}\mu(A,U) z^{r\dim(A)}. \\
   & = \sum_{V \leq E} p(M.V;z) \sum_{j=0}^{\dim(U\cap V^\perp)} \qbin{\dim(U\cap V^\perp)}{j}{q}(-1)^{\dim(U)-j} q^{\binom{\dim(U)-j}{2}}z^{jr}. \qedhere
\end{align*}
\end{proof}

We now have the following MacWilliams identity, relating the weight enumerators of $M$ and $M^*$. This version of the identity, or rather its corollary, will be used to prove the main theorem of Section~\ref{sec:designs}.
\begin{Theorem}\label{th:macw}
	Let $s \in \{0,\ldots ,n\}$. Then
	$$\sum_{i=0}^{n-s} \qbin{n-i}{s}{q} A_M(i;z) = z^{{\rho(E)}-rs} \sum_{i=0}^s {\qbin{n-i}{s-i}{q}} A_{M^*}(i;z). $$
\end{Theorem}

\begin{proof}
	We start with the left-hand-side of the equation and rewrite it, noting that
	$\qbin{n-i}{s}{q}=\qbin{n-i}{n-s-i}{q}$ counts the number of $(n-s)$-dimensional subspaces of $E$
	that contain a fixed space of dimension $i$.
	This yields:
	\begin{eqnarray*}
	   \sum_{i=0}^{n-s} \qbin{n-i}{s}{q} A_M(i;z) &=& 
	   \sum_{i=0}^{n-s} \qbin{n-i}{n-s-i}{q} \sum_{X:\dim(X)=i} p(M.X;z)\\
	    & = & \sum_{{U\,:\,} \dim(U)=n-s} \sum_{X \leq U} p(M.X;z).\\
    \end{eqnarray*}  
    From Lemma \ref{lem:charMM^*}, this gives:
\begin{eqnarray*}
	\sum_{i=0}^{n-s} \qbin{n-i}{s}{q} A_M(i;z)   & = &
	   \sum_{{U\,:\, }\dim(U)=n-s} z^{{\rho(E)}-r\dim(U^\perp)} \sum_{X \leq U^\perp} p(M^*.X;z) \\
	   & = & \sum_{{V\,:\,} \dim(V)=s} z^{{\rho(E)}-rs} \sum_{X \leq V} p(M^*.X;z)\\
	    & = & z^{{\rho(E)}-rs} \sum_{i=0}^s {\qbin{n-i}{s-i}{q}} \sum_{X \leq {E\,:\,} \dim(X)=i} p(M^*.X;z)\\
	   &=& z^{{\rho(E)}-rs} \sum_{i=0}^s {\qbin{n-i}{s-i}{q}} A_{M^*}(i;z).
	\end{eqnarray*}
\end{proof}

{
Theorem \ref{th:macw} shows that the weight enumerator of a $q$-polymatroid and that of its dual are related by invertible $q$-Pascal matrices. The minors of such matrices have been studied as $q$-analogues of the classical Pascal matrices. We will use the following result from \cite[Theorem 2.2]{zhizheng}.
\begin{Lemma}\label{lem:qpascal}
Let $r_1,\ldots,r_n$ be a sequence of non-negative integers. We have
\[
\det\left(\qbin{r_i}{j-1}{q}\right)_{1\leq i,j \leq n} = q^{\binom{n}{2}} \prod_{1 \leq i < j \leq n} \frac{q^{r_j}-q^{r_i}}{q^j-q^i}.
\]
\end{Lemma}
The next corollary (c.f. \cite[Corollary 3.2]{BRS2009}) is the main device used to prove Theorem \ref{th:weightedtdesign}, which identifies sufficiency criteria for the existence of weighted subspace designs arising from the dependent spaces of a $q$-polymatroid (c.f. \cite[Theorem 3.3 ]{BRS2009}).
We remark that the reasoning used here is similar to that of the original Assmus-Mattson Theorem and its generalizations. For any positive integer $\ell$, we write $[\ell]:=\{1,\ldots,\ell\}$.}

\begin{Corollary}\label{cor:macw}
	Let $S \subseteq \{1,\ldots ,n\}$.
	The pair of lists
	$$[A_{M^*}(i;z): |S| \leq i \leq n] \text{ and } [A_M(j;z):j \in S],$$
	is determined uniquely by the pair of lists
	$$[A_{M^*}(i;z): 1\leq i \leq |S|-1] \text{ and } [A_M(j;z):j \in [n] - S].$$
\end{Corollary}

\begin{proof}
    {Let $A_M(z):=\left( A_M(i;z) \right)_{0\leq i \leq n}$ and let 
    $A_{M^*}(z):=\left( A_{M^*}(i;z) \right)_{0\leq i \leq n}$. Note that $A_M(0;z)=A_{M^*}(0;z)=1$, and in particular are known.}
	From Theorem \ref{th:macw}, we have the matrix equation
	$$\left(\qbin{n-i}{s}{q}\right)_{{0 \leq i,s \leq n}} A_M(z) 
	= {\rm diag}(z^{{\rho(E)}-rs})_{0 \leq s \leq n}\left({\qbin{n-i}{n-s}{q}}\right)_{{0 \leq i,s \leq n}} A_{M^*}(z).$$
	{Let $t=|S|$ and write  $S=\{\ell_1,\ldots ,\ell_t\}$. 
	By Lemma \ref{lem:qpascal}, we have
	\[
	   \det\left(\qbin{n-\ell_i}{s-1}{q}\right)_{1 \leq i,s \leq t} = q^{\binom{t}{2}} \prod_{1 \leq i < s \leq t} \frac{q^{n-\ell_s}-q^{n-\ell_i}}{q^s-q^i},
	\]
	which is non-zero, as the $\ell_i$ are distinct. 
	}
	{
	Now suppose that the coefficients $A_M(j;z)$ are known for $j \notin S$ and that the $A_{M^*}(j;z)$ are known for $0 \leq j \leq t-1 $. Then we can solve for the unknown $A_M(j;z)$ via
	\begin{small}
	\begin{align*}
	    \left( A_M(\ell_i;z) \right)_{1 \leq i \leq t} 
	    &= \left(\qbin{n-\ell_j}{s-1}{q}\right)_{1 \leq j,s \leq t}^{-1} \left( {\rm diag}(z^{{\rho(E)}-rs})_{0 \leq s \leq t-1}\left({\qbin{n-j}{n-s}{q}}\right)_{{0 \leq s,  j \leq t-1}} \left( A_{M^*}(j;z) \right)_{0\leq j \leq t-1} \right. \\
	    &-\left. \left(\qbin{n-j}{s-1}{q}\right)_{
	    \substack{1 \leq s \leq t\\j \in \{0,\ldots, n\} - S }} \left( A_M(j;z) \right)_{j \in \{0,\ldots, n\} - S} \right). 
    \end{align*}
	\end{small}
	Once the list $[A_M(j;z):j \in S]$ is determined, since $A_M(z)$ is now known, Theorem \ref{th:macw} can be applied to retrieve $[A_M^*(i; z): t \leq i \leq n]$.
	}
\end{proof}

\section{Weighted Subspace Designs from {\em q}-Polymatroids}\label{sec:designs}

\subsection{Weighted Subspace Designs}

In \cite{BRS2009}, the authors give a definition of a {\em weighted subspace design}, which generalizes a $t$-design. A $t$-$(n,k,\lambda)$ design, with $t,k,\lambda$ positive integers, is a collection of $k$-subsets of an $n$-set (called blocks) with the property that every $t$-subset of the $n$-set is contained in exactly $\lambda$ blocks. A $q$-analogue of this notion is that of a $t$-design over $\FF_q$, which is a collection of $k$-dimensional subspaces of $E$ called blocks, with the property that every $t$-dimensional subspace of $E$ is contained in the same number of blocks. Similarly, there is a $q$-analogue of a weighted $t$-design.  

\begin{Definition}\label{weighted-t-design}
Let ${\mathbb G}$ be an additive group, $t,k$ positive integers and $\lambda \in \mathbb{G}$. A weighted $t$-$(n,k,\lambda; q)$ design $\mathcal{D}$ is a triple  $(E,\mathcal{B},f)$ for which $\mathcal{B}$ is a collection of $k$-subspaces  of $E$ (called blocks) and $f:\mathcal{B} \mapsto \mathbb{G}$ is a weight function such that for all $t$-spaces $T$ of $E$, $\displaystyle \sum_{B:\\ T\leq B} f(B)= \lambda$. We say that $\D$ is a {\bf weighted subspace design} or is a {\bf weighted design over} $\FF_q$.
\end{Definition}

A subspace design (a design over $\FF_q$) can be interpreted as a weighted subspace design with the weight function $f(B):=1$ for all $B\in \mathcal{B}$, and ${\mathbb G} = ({\mathbb Z},+)$. For an excellent survey on subspace designs, see \cite{bks18}. In general, obtaining new subspace designs is a difficult problem, often highly dependent on computer search, which is exacerbated by the number of subspaces involved (which is exponential in comparison to classical designs for the same parameters). For example, it is not yet known if a $3$-$(8,4,1;2)$ subspace design exists; such a design would have $6477$ blocks, chosen from an ambient space having $200{,}787$ $4$-dimensional subspaces. Its classical analogue, the extended Fano plane, has $14$ blocks, chosen from a collection of $70$ $4$-sets.  In \cite{WINEpaper1}, a construction of a $q$-analogue of a perfect matroid design ($q$-PMD) was given, which is a $q$-matroid for which all flats of the same dimension have the same rank. This $q$-PMD yields a construction of a subspace design from a $q$-Steiner system. In the following sections we will show another way that subspace designs and weighted subspace designs can arise from $q$-polymatroids satisfying certain rigidity properties. \\

The \textbf{intersection numbers} of a weighted subspace design are important invariants and can be used to establish non-existence results. Their values are the same as for subspace designs; see, for example \cite[Fact 1.5]{KP} or \cite{suzuki1990inequalities}.

\begin{Theorem}\label{intersection-numbers}
  Let  $(E,\mathcal{B},f)$ be a $t$-$(n,k,\lambda; q)$ weighted subspace design and let
	 $I, J$ be two subspaces of $E$ of dimension $i$ and $j$, respectively, such that $I \cap J=\{0\}$. If $i+j \leq t$, then 
	\begin{align*}
	\sum_{B\in { \mathcal{B}\,:\,} I \leq B, B \cap J=\{0\} }f(B) = q^{(k-i)j} \qbin{n-i-j}{k-i}{q}\qbin{n-t}{k-t}{q}^{-1} \lambda.
	\end{align*}
	In particular, this number is independent of the choice of $I$ of dimension $i$ and $J$ of dimension $j$. We denote it by $\lambda_{i,j}$.
\end{Theorem}
\begin{proof}
 If $X$ is a subspace of $E$ of dimension $x \leq t$, then since $(E,\mathcal{B},f)$ is a weighted subspace design, we have
\begin{align}\label{first-eqn-of-2.6}
	\qbin{k-x}{t-x}{q} \sum_{B\in {\mathcal{B}\,:\, X}\leq B} f(B) &= \sum_{B\in { \mathcal{B}\,:\, X}\leq B} \quad
	\sum_{{T\,:\, X}\leq T\leq B, \dim(T)=t} f(B)  ,\\ &=
	\sum_{{T\,:\, X}\leq T\leq B, \dim(T)=t} \quad \sum_{B\in {\mathcal{B}\,:\, T}\leq B} f(B) \nonumber =
	\qbin{n-x}{t-x}{q} \lambda .
\end{align}
Now restrict to a subspace $X$ of the form $X=I+L$ for some $L \leq J$ of dimension $s$. Then $I \cap L = \{0\}$ and 
$\dim(I+L) = i+s$ and so (\ref{first-eqn-of-2.6}) becomes:
\begin{align*}
   g(L):=\sum_{B\in {\mathcal{B}\,:\, I+L}\leq B} f(B)  =
	\qbin{n-(i+s)}{t-(i+s)}{q}\qbin{k-(i+s)}{t-(i+s)}{q}^{-1}\lambda. 
\end{align*}
Define $\displaystyle h(K) = \sum_{B \in {\kB \,:\, I} \leq B , B \cap J = K} f(B)$, for each $K \leq J$. Then we have that
\[g(L)=\sum_{K\; :\;  L \leq K \leq J} h(K) \]
and so, by M\"{o}bius inversion on the lattice $\kL(J)$,
\[ h(L)=\sum_{K\; :\; L\leq K\leq J} \mu(L,K) g(K).\]
Substituting $L=\{0\}$ now gives 
\begin{align*}
   \sum_{B\in { \mathcal{B}\,:\,} I \leq B, B \cap J=\{0\} }f(B) & = h(\{0\}) = \sum_{K \leq J} \mu(0,K) g(K) \\
   & = \sum_{s=0}^j \qbin{j}{s}{q} (-1)^s q^{\binom{s}{2}}\qbin{n-i-s}{t-i-s}{q}\qbin{k-i-s}{t-i-s}{q}^{-1}\lambda \\
   & = \lambda \qbin{n-t}{k-t}{q}^{-1} \sum_{s=0}^j (-1)^s q^{\binom{s}{2}}  \qbin{j}{s}{q} \qbin{n-i-s}{k-i-s}{q}\\
   & = \lambda \qbin{n-t}{k-t}{q}^{-1} q^{j(k-i)}\qbin{n-i-j}{k-i}{q}
\end{align*}
The third line follows from applying Equation $\eqref{eq:qbinid}$ with 
$a=n-i-s$, $b= k-i-s$, $c=t-i-s$, 
and the last equality follows from Lemma \ref{counting_subspace}.
   
\end{proof}

The proof outlined above is a direct $q$-analogue of \cite[Theorem 2.6]{BRS2009}.
The intersection numbers for subspace designs were given in \cite{Cameron,suzuki1990inequalities}, for which the authors proposed an inductive argument.

We have the following constructions of weighted subspace designs from a given one (c.f. \cite{KP,suzuki1990inequalities}).

\begin{Corollary}\label{cor:weighteddesign}
	Let $\mathcal{D}:=(E, \mathcal{B},f)$ be a weighted $t$-$(n,k,\lambda;q)$ design. 
	\begin{enumerate}
		\item  For $0\leq i\leq t$, $\D$ is an
		$i$-$(n,k,\lambda_i; q)$ weighted subspace design with
		\[
		\lambda_i=\qbin{n-i}{k-i}{q}\qbin{n-t}{k-t}{q}^{-1}\lambda ~
		{=\qbin{n-i}{t-i}{q}\qbin{k-i}{t-i}{q}^{-1}\lambda}.
		\]
		\item  Define $\mathcal{B}^{\bot}:=\{B^{\bot}: B\in \mathcal{B}\}$. 
		{If $k \leq n-t$}
		then $\mathcal{D}^{\bot}=(E, \mathcal{B}^{\bot},{f^\perp})$ is a $t$-$(n,n-k,\lambda^\perp; q)$
		weighted subspace design with ${f^\perp(B^\perp):=f(B)}$ for all $B \in \mathcal{B}$ and 
		$\lambda^\perp:=\qbin{n-k}{t}{q}\qbin{k}{t}{q}^{-1}\lambda$.
	\end{enumerate}
\end{Corollary}
\begin{proof}
    To see that (1) holds, apply Theorem \ref{intersection-numbers} with $\lambda_i: = \lambda_{i,0}$. Let $I$ be an $i$-dimensional subspace of $E$.
    We have $\lambda_{i,0} = \sum_{B \in {\kB\,:\, I} \leq B} f(B)=
    \qbin{n-i}{k-i}{q}\qbin{n-t}{k-t}{q}^{-1}$. The rest follows from {Equation} (\ref{eq:qbinid}).
  
    {We will compute the value $\lambda^\perp$. 
    A $t$-dimensional subspace $T$ is contained in $B^\perp \in \kB^\perp$ if and only if $B \leq T^\perp$.
    Now consider the set $S :=\{(B,X): B \in \mB, \dim(X) = n-t, B \leq X \}$.
    We will compute the sum of the $f(B)$ over all pairs $(B,X)$ in $S$ in two ways.
    On the one hand, we have:
    \begin{align*}
        \sum_{(B,X) \in S} f(B) &= 
        \sum_{B \in \mB}  \sum_{\substack{X: B \leq X \\ \dim(X)=n-t}} f(B) =
        \qbin{n-k}{n-t-k}{q} \sum_{B \in \mB}  f(B)=
        \qbin{n-k}{t}{q} \lambda_{0,0}\\
        &=\qbin{n-k}{t}{q} \qbin{n}{k}{q}\qbin{n-t}{k-t}{q}^{-1} \lambda =
        \qbin{n-k}{t}{q} \qbin{n}{t}{q}\qbin{k}{t}{q}^{-1} \lambda.
    \end{align*}    
The last equality follows from applying Equation \eqref{eq:qbinid} with $a=n$, $b=k$, $c=t$. On the other hand,
\[ \sum_{(B,X) \in S} f(B) =\sum_{\substack{X \leq E: \\ \dim(X)=n-t}} \sum_{B \in \mB, B \leq X} f(B)= \qbin{n}{t}{q} \sum_{B \in \mB, B \leq X} f(B). \]
    It follows, by comparing the two right-hand sides, that 
    \[\displaystyle \lambda^\perp:=\sum_{B^\perp \in \mB^\perp, T \leq B^\perp} f(B) 
    =\sum_{B \in \mB, B \leq T^\perp} f(B)  = \qbin{n-k}{t}{q} \qbin{k}{t}{q}^{-1} \lambda. 
    \]}
    \end{proof}

\subsection{Subspace Designs from {\em q}-Polymatroids}

We now present criteria for the existence of a weighted subspace design arising from the dependent spaces of a $q$-polymatroid. The approach is in essence a generalization of the original argument given by Assmus and Mattson \cite{AM69}. To do this, we obtain a $q$-analogue of \cite[Theorem 3.3]{BRS2009}.
Throughout this section we let $\FF$ denote an arbitrary field. 
Since $p(M;z) \in \ZZ[z]$, it gives a well-defined function on any field, viewed as a $\ZZ$-module. 
We define the following (c.f. \cite{BRS2009}).

\begin{Definition}\label{def:DRd_M}
Let $\theta \in \FF$. We define:
\begin{itemize}
	\item $D_M(i;\theta):=\{ X \leq E:\dim(X)=i, p(M.X;\theta) \neq 0\}$,
	\item $R_M(t;\theta):=\{ j \in \{1,\ldots,n-t\}: A_{M^*}(j;\theta) \neq 0\}$,
	\item $d_M:=\min\{ \dim(X):X \leq E, X \text{ is a cocircuit of }M \}$.
\end{itemize}
\end{Definition}

The sets $D_M(i;\theta)$ will, in certain circumstances, form the blocks of weighted subspace designs.

\begin{Proposition}\label{prop:D_M(d_M;z)}
	Let $\theta \in\FF$ such that $\theta^{s} \neq 1$ for any $s \in \{1,\ldots, r\}$.
	{For each $i \in \{1,\ldots,n\}$
	every member of $D_M(i;\theta)$ is a dependent space of $M^*$.} Moreover $D_M(d_M;\theta)$ is precisely the set of circuits of $M^*$ of dimension $d_M$.
\end{Proposition}

\begin{proof}
If $A \in D_M(i;\theta)$, then $p(M.A;\theta) \neq 0$, which by Lemma \ref{lem:pm.tneq0} means that $A$ is a dependent space of $M^*$.
We show that for circuits of $M^*$ (i.e., minimal dependent spaces of $M^*$) the converse also holds.
By Lemma \ref{lem:pm.tcirc}, for any circuit $X$ of $M^*$ we have $p(M.X;z)=z^{\ell(X^\perp)}-1$. Since $X$ is not independent in $M^*$, $\ell(X^\perp)>0$ by Lemma \ref{lem:indlperp0}.
Therefore, by our choice of $\theta$, we have that $p(M.X;\theta)=\theta^{\ell(X^\perp)}-1\neq0$ and so $X \in D_M(\dim(X);\theta)$. In particular, $D_M(d_M;\theta)$ is precisely the set of all circuits of $M^*$ of dimension $d_M$. 
\end{proof}

We will now present the main results of this section: Theorem \ref{th:weightedtdesign} and its two corollaries. Together they form a $q$-analogue of \cite[Theorem 3.3]{BRS2009}.

\begin{Theorem}\label{th:weightedtdesign}
	Let $\theta \in\FF$ such that $\theta^{s} \neq 1$ for any $s \in \{1,\ldots, r\}$. Let $t<d_M$ be a positive integer. Suppose that $\sigma^*:=|R_M(t;\theta)| \leq d_M - t$ and suppose further that
	for each $t$-dimensional subspace $T$ and $j\leq n-t$, if   
	$A_{M^*}(j;\theta) = 0 $, then  $A_{M^*/T}(j;\theta) = 0.$
	Then $(E,D_M(d_M;\theta),f)$ is a weighted $t$-design over $\mathbb{F}_q$
	with $f(X):=p(M.X;\theta)$. 
\end{Theorem}

\begin{proof}
	Let $T$ be a $t$-dimensional subspace of $E$. Since $t<d_M$, $T$ is independent in $M^*$.
	By Lemma \ref{lem:MtoM/T}, any dependent space $A$ of $M^*/T$ has the form $A=B/T$ for a dependent space $B$ of $M^*$.
	Therefore, {for any such $A$ and $B$ we have} 
	\begin{eqnarray}\label{eq:dMt}
		\sigma^* \leq  d_M-t \leq \dim(B)-t = \dim(A).
	\end{eqnarray}
	In other words, no dependent space of $M^*/T$ has dimension less than $\sigma^*$.
	By Lemma \ref{lem:pm.tneq0}, if $X$ is {non-trivial} and independent in 
	$M^*/T${, then} $p((M^*/T)^*.X;\theta)=0$. 
	Therefore, 
	\begin{align*}\displaystyle A_{(M^*/T)^*}(i;\theta)=\sum_{X \leq {E/T\,:\, \dim(X)=i}} p((M^*/T)^*.X;\theta)=0, \text{ for all }{1}\leq i \leq  \sigma^*-1.
	\end{align*}
	By hypothesis, $A_{M^*/T}(j;\theta) =0$
	for all $j \notin R_M(t;\theta)$, {$j\leq n-t$} and 
	so the coefficients,
	\begin{align*} [A_{M^*/T}(j;\theta) :  j \notin R_M(t;\theta), {j\leq n-t}]\text{ and } [A_{(M^*/T)^*}(i;\theta): {1}\leq i \leq  \sigma^*-1],
	\end{align*}
	are known.
	Now write $S = R_M(t;\theta)$ and apply Corollary \ref{cor:macw} to see that the coefficients 
	\begin{align*} [A_{M^*/T}(j;\theta) : j \in R_M(t;\theta)] \text{ and } [A_{(M^*/T)^*}(i;\theta): \sigma^* \leq i \leq n-t] 
	\end{align*}
	are uniquely determined and independent of our choice of $T$ of dimension $t$.
	It follows that the $A_{(M^*/T)^*}(i;\theta)$ 
	are uniquely determined for $0\leq i \leq n-t$. 
	We will now show that
	\begin{align*}
		\sum_{X \in {D_M(d_M;\theta)\,:\, T} \leq X} p(M.X;\theta) =  A_{(M^*/T)^*}(d_M-t;\theta),  
	\end{align*} 
	which will establish that $(E,D_M(d_M;\theta),f)$ is a weighted $t$-design over $\mathbb{F}_q$
	with $f(X):=p(M.X;\theta)$.
	
	{
	We claim there is a one-to-one correspondence between the members of $D_M(d_M;\theta)$ that contain $T$ and the members of $D_{(M^*/T)^*}(d_M-t;\theta)$.}
	Let $B$ be a circuit of $M^*$ that contains $T$ such that $\dim(B) = d_{M}$. From Lemma \ref{lem:MtoM/T},
	$B/T$ is a circuit of $M^*/T$ and $\dim(B/T) = \dim(B)-t = d_M - t$. Conversely, if $A$ is a circuit
	of $M^*/T$ satisfying $\dim(A) = d_M - t$, then $A = B/T$ for a dependent space $B$ of $M^*$ of dimension $\dim(B)=d_M$, which is therefore a circuit of $M^*$, as it has minimal dimension. 
	By Proposition \ref{prop:D_M(d_M;z)}, $D_M(d_M;\theta)$ is the set of all cocircuits of $M$ of dimension $d_M$ and hence there is a one-to-one correspondence between the members of $D_M(d_M;\theta)$ that contain $T$ and the circuits of $M^*/T$ of dimension $d_M-t$.
	By Equation~\eqref{eq:dMt}, any dependent space of $M^*/T$ of dimension $d_M-t$ is a circuit of $M^*/T$ and hence is a member of $D_{(M^*/T)^*}(d_M-t;\theta)$. 
	This establishes the claim.
	 
	From Corollary \ref{cor:carpolyMTperp}, for any circuit $X/T$ of $M^*/T$ we have 
	$$p((M^*/T)^*.(X/T);\theta)
	= \theta^{\ell(X^\perp)}-1.$$
	 Therefore,  
	\begin{eqnarray*}
	\sum_{X \in D_M(d_M;\theta): T \leq X} p(M.X;\theta) 
	&=& \sum_{X \in D_M(d_M;\theta)\,:\, T \leq X} (\theta^{\ell(X^\perp)}-1), \\
	&=&\sum_{{X/T \in D_{(M^*/T)^*}(d_M-t;\theta)}} (\theta^{\ell(X^\perp)}-1),\\	
	&=&\sum_{X/T \leq{ E/T\,:\, } \dim(X/T)=d_M-t}p((M^*/T)^*.(X/T);\theta),\\ 
	&=&A_{(M^*/T)^*}(d_M-t;\theta),
	\end{eqnarray*}
	which is independent of our choice of $T$ of dimension $t$.
	It follows that $(E,D_M(d_M;\theta),f)$ is a weighted $t$-design over $\mathbb{F}_q$
	with $f(X):=p(M.X;\theta)$. 
\end{proof}

\begin{Remark}\label{rem:Aisindep}
	In the proof of Theorem \ref{th:weightedtdesign}, we saw that with the hypothesis of the theorem, that the $A_{(M^*/T)^*}(i;\theta)$ (and therefore the $A_{M^*/T}(i;\theta)$)
	are uniquely determined for $0\leq i \leq n-t$.  
	By Lemma \ref{lem:dualMTperp}, it follows that the $A_{{M|T^\perp}}(i;\theta)$ are uniquely determined for $0\leq i \leq n-t$.	
\end{Remark}
	
\begin{Corollary}\label{firstcorollary}
Let $\theta \in\FF$ such that $\theta^{s} \neq 1$ for any $s \in \{1,\ldots, r\}$. Let $t<d_M$ be a positive integer. Suppose that $\sigma^*:=|R_M(t;\theta)| \leq d_M - t$ and suppose further that
for each $t$-dimensional subspace $T$ and $j\in \{d_M,\ldots, n-t\}$, if
$A_{M^*}(j;\theta) = 0$ then $A_{M^*/T}(j;\theta) = 0.$
Then for each $j \in \{d_M,\ldots ,n-t\}$, $(E,D_M(j;\theta),f)$ is a weighted $t$-design over $\mathbb{F}_q$
with $f(X):=p(M.X;\theta)$.
\end{Corollary}
\begin{proof}
    We will prove by induction on $w \in \{d_M,\ldots,n-t\}$ that $(E,D_M(w,\theta),f)$ is a weighted $t$-design. The first step was proved in Theorem~\ref{th:weightedtdesign}. 
 Suppose now that $(E,D_M(j,\theta),f)$ is a weighted $t$-design
 for each $j \in \{d_M,\ldots, w-1\}$. We will show that $(E,D_M(w,\theta),f)$ is also a  weighted $t$-design.
 
 Let $T\leq E$ have dimension $t$.
 We will show that the following sum depends only on $t$:
 $$\sum_{W\in D_M(w;\theta),\: T \leq W} p(M.W;\theta) = \sum_{W\, :\, T \leq W,\; \dim(W)=w} p(M.W;\theta).$$ 
 Note that since $D_M(w;\theta)$ is the set of $w$-dimensional subspaces of $E$ for which $p(M.W;\theta)\neq 0$, the above equality holds.
 From Lemma \ref{lem:2.11}, for any $T \leq W \leq E$ we have that 
 \begin{equation*}
 	p({M|T^\perp}/W^\perp;\theta) = \sum_{A\,:\, A+T =W} p(M.A;\theta).
\end{equation*}

Let $\phi: {\mathcal L}(E/T) \longrightarrow {\mathcal L}(T^\perp)$ be defined by $\phi(A/T)= (A^\perp)^{\perp(T^\perp)}$, for each $A \leq E$ such that $T \leq A$. (This map was already used in Lemma \ref{lem:dualMTperp}.) Let $X$ be such that $\phi(W/T)=X$.
Then $(W^\perp)^{\perp(T^\perp)}=X$ and so $X^{\perp(T^\perp)} = W^\perp$. Therefore,
	${M|T^\perp}.X = (M|T^\perp)/W^\perp$.
	Clearly, $\dim(W) = \dim(X)+\dim(T)$ and so if $T$ is a $t$-dimensional space{, then}:
	\[A_{{M|T^\perp}}(j;z) = \sum_{X \leq {T^\perp\,:\, }\dim(X)=j} p({M|T^\perp}.X;z) = \sum_{W \leq { E\,:\,} T \leq W,\; \dim(W)=j+t} p({M|T^\perp}/W^\perp;z).
	\]

Therefore we have
 \begin{equation*}
 	A_{{M|T^\perp}}(w-t;\theta)= \sum_{\substack{W\,:\,T\leq W,\\ \dim(W)=w}} p({M|T^\perp}/W^\perp;\theta)
 	 =  \sum_{\substack{W\,:\,T\leq W,\\ \dim(W)=w}} \:\sum_{{A\,:\, A}+T =W} p(M.A;\theta)
 \end{equation*}
 For any $I\leq T$, we write $I_T$ to denote an arbitrary fixed subspace of $T$ satisfying $I\oplus I_T=T$.
 Clearly, if $I = A \cap T$ we have $I\leq A$ and $A \cap I_T=\{0\}$. Conversely, if $I \leq A$ and $I_T \cap A = \{0\}$ then $A\cap T =A\cap(I+I_T)=I$.
 Moreover, if $W$ is a $w$-dimensional subspace for which $A+T=W$ and $A \cap T= I$, then $\dim(A)=w-t+\dim(I)$.
 Therefore, we can rewrite the double summation as follows:
 \begin{eqnarray*}
     A_{{M|T^\perp}}(w-t;\theta)&=& 
     \sum_{i=0}^t \:
     \sum_{\substack{I\, :\, I \leq T,\\ \dim(I)=i}} \:
     \sum_{\substack{A\, :\, I \leq A,\: I_T \cap A =\{0\},\\ \dim(A)=w-t+i} } p(M.A;\theta)\\
     & = & \sum_{\substack{A\,:\, T \leq A,\\ \dim(A)=w}} p(M.A;\theta) +  
     \sum_{i=0}^{t-1} \:\sum_{\substack{I\,:\,I \leq T \\ \dim(I)=i}} 
     \:\sum_{\substack{A\, :\, I \leq A,\: I_T \cap A =\{0\},\\ \dim(A)=w-t+i} } p(M.A;\theta).
 \end{eqnarray*}
 Let $I\leq T$ such that $\dim(I)=i<t$ and so $d_M-t\leq w-t\leq w-t+i\leq w-1$.
 By hypothesis, for each $1\leq j \leq w-1$, $(E,D_M(j,\theta),f)$ is a weighted $t$-design with $f(X):=p(M.X;\theta)$, and so by
 Theorem \ref{intersection-numbers},
 $$\sum_{\substack{A\, :\, I \leq A,\: I_T \cap A =\{0\},\\ \dim(A)=w-t+i} } p(M.A;\theta) = \Lambda^w_{i,t-i}(M;\theta),$$ 
 for $\Lambda^w_{i,t-i}(M;\theta)$ that depend only on $t,w,i$. 
 It follows that 
 \begin{eqnarray*}
 	 \sum_{\substack{W\,:\,T \leq W,\\ \dim(W)=w}} p(M.W;\theta)
 	  &=& 
 	      A_{{M|T^\perp}}(w-t;\theta)-
 	      \sum_{i=0}^{t-1}\qbin{t}{i}{q} \Lambda^w_{i,t-i}(M;\theta).
 \end{eqnarray*}
  By Remark \ref{rem:Aisindep},
  $A_{{M|T^\perp}}(w-t;\theta)$ is independent of our choice of $T$ of dimension $t$ and so the result follows.
\end{proof}

\begin{Corollary}\label{secondcorollary}
Let $\theta \in\FF$ such that $\theta^{{s}} \neq 1$. Let $t<d_M$ be a positive integer. Suppose that $\sigma^*:=|R_M(t;\theta)| \leq d_M - t$ and suppose further that
for each $t$-dimensional subspace $T$ and $j\leq n-t$, if $A_{M^*}(j;\theta) = 0$, then $A_{M^*/T}(j;\theta) = 0.$
Then for each $j \in \{d_{M^*},\ldots,n-t\}$, $(E,D_{M^*}(j;\theta),f^*)$ is a weighted $t$-design over $\mathbb{F}_q$
with $f^*(X):=p(M^*.X;\theta)$ for all subspaces $X \leq E$.
\end{Corollary}

\begin{proof}
For each $d_{M^*}\leq j\leq n-t$, define the set $\mathcal{D}_j:=\{X^\perp: X\in D_{M^*}(j;\theta)\}$. 
Let $T$ be a $t$-dimensional subspace of $E$.
Then for each $j$ we have:
\begin{eqnarray*}
\sum _{X \in D_{M^*}(j,\theta): T\leq X^{\perp}} p(M^*.X;\theta)
=\sum _{{X\,:\,X} \leq T^\perp, \dim(X)=j}p(M^*.X;\theta) = 
\sum _{{X\,:\,X} \leq T^\perp, \dim(X)=j}p(M^*/X^\perp;\theta).
\end{eqnarray*}
Now for each $X \leq T^\perp$ we have $(E/T) \big/ (X^\perp/T) \cong E/X^\perp$ and it is easy to see that the corresponding $q$-polymatroids are lattice-equivalent. 
Then, applying Lemma \ref{lem:dualMTperp} and using the fact that $\phi(X^\perp /T)^{\perp(T^\perp)} = X$,
we get
\[
M^*.X  \cong (M^*/T)\big/ (X^\perp/T)  \cong ({M|T^\perp})^* / \phi(X^\perp/T) 
= ({M|T^\perp})^*.\phi(X^\perp/T)^{\perp(T^\perp)}
\cong ({M|T^\perp})^*.X.
\] 
Therefore,
\[
 \sum _{X \in \D_j: T\leq X} p(M^*.X;\theta) =\sum _{{X\,:\,T} \leq X, \dim(X)=n-j}p(({M|T^\perp})^*.X;\theta) = A_{({M|T^\perp})^*}(n-j;\theta). 
\]
From Remark \ref{rem:Aisindep}, $A_{({M|T^\perp})^*}(n-j;\theta)$ is independent of the choice of $T$ of dimension $t$. It follows that $(E,\D_j,f^*)$ is a weighted subspace design with $f^*$ defined by $f^*(X)=p(M^*.X;\theta)$ for each $X \leq E$. The result now follows by Corollary \ref{cor:weighteddesign}: the required subspace design is the dual of $(E,\D_j,f^*)$.
\end{proof}

\begin{Remark}
   The 
   results of Proposition \ref{prop:D_M(d_M;z)}, Theorem \ref{th:weightedtdesign} and Corollaries \ref{firstcorollary} and \ref{secondcorollary} all hold with indeterminate $z$ in place of a specific choice of $\theta$ in $\FF$.
   In particular, $p(M.X;z)$ is a non-zero polynomial in $\ZZ[x]$ for any cocircuit $X$ of $M$.
\end{Remark}

In general, a $(q,r)$-polymatroid $M$ may satisfy the hypothesis of Corollary \ref{firstcorollary} for one choice of $\theta$, but fail for another choice. However, if the hypothesis holds for indeterminate $z$, then a weighted $t$-design over $\FF_q$ can be constructed for any choice of $\theta$ that doesn't vanish on $p(M.X;z)$ for a cocircuit $X$ of $M$.

\begin{Example}
	Let $M=(E,\rho)$ be the uniform $q$-matroid $U_{k,n}$, as described in Example \ref{exunif}. We will show that this $q$-matroid satisfies the hypothesis of Corollary \ref{firstcorollary} with indeterminate $z$ in place of a specific choice of $\theta$ in some field $\FF$.
	
	The dual $q$-matroid $M^*=(E,\rho^*)$ is the uniform $q$-matroid $U_{n-k,n}$, whose independent spaces are exactly those of dimension $n-k$ or less, and for which all other spaces are dependent and have rank $n-k$. 
	Therefore, every cocircuit of $M$ has dimension $d_M = n-k+1$.
	Now $p(M^*.X;z)=0$ for all subspaces $X$ such that $1 \leq \dim(X) \leq k$, as these are the independent spaces of $M$ (see Lemma \ref{lem:pm.tneq0}), and so 
	$A_{M^*}(i;z) = 0$ for all $i \in \{1,\ldots ,k\}$.
	Therefore for any $t\leq d_M-1=n-k$, we have $R_M(t;z) \subseteq \{k+1,\ldots ,n-t\}$ and so $|R_M(t;z)| \leq n-t-k \leq d_M -t$.
	
	We now show that for any indeterminate $z$ the condition holds that for all $j\leq n-t$, if $A_{M^*}(j;z) = 0$, then $A_{M^*/T}(j;z) = 0$. 
	Let $T$ be a $t$-dimensional subspace of $E$ and recall $1 \leq t \leq n-k$. By Lemma \ref{lem:M/T} 
	we have
	\[ 
	A_{M^*/T}(j;z)= \sum_{X \leq {T^\perp\,:\, \dim(X)=j}} p(M^*.X;z).
	\]
	Since $p(M^*.X;z)=0$ for all subspaces $X$ such that $1 \leq \dim(X) \leq k$, we have that $A_{M^*}(j;z) = 0$ for all $j\in\{1,\ldots,k\}$.
	
	Next, we consider the case $k+1\leq j\leq n-t$. Let $X \leq E$ be a subspace of dimension at least $k+1$. We claim that the $q$-matroid $M^*.X$ has no loops, in which case by Lemma \ref{lem:monic},  
	$p(M^*.X;z)$ will be a monic polynomial of degree $n-k-\dim(X^\perp)$ and hence $A_{M^*}(j;z) \neq 0$ for $k+1 \leq j \leq n-1$, i.e. the condition holds vacuously. Consider a subspace $U$ that strictly contains $X^\perp$.
	Since $M^*=U_{n-k,n}$, we have $\rho^*(X^\perp)=\dim(X^\perp)$ and so
	\begin{align*}
	\rho_{M^*/X^\perp}(U/X^\perp) & =\rho^*(U)-\rho^*(X^\perp) \\
	& = \rho^*(U)-\dim(X^\perp) \\
	& = \dim U-\rho(E)+\rho(U^\perp)-\dim( X^\perp) \\
	& = \dim U-\dim X^\perp+\rho(U^\perp)-k
	\end{align*}
	We have that $\rho(U^\perp)=\min\{\dim U^\perp,k\}$. Substituting both cases in the equation above and using that $\dim X-k\geq 1$ and $\dim U-\dim X^\perp\geq 1$, respectively, we find that $\rho_{M^*/X^\perp}(U/X^\perp)\geq 1$. This implies that the $q$-matroid $M^*.X$ has no loops.

	We conclude that $M=U_{k,n}$ satisfies the hypothesis of Corollary \ref{firstcorollary} for indeterminate $z$. Therefore 
	for any $\theta \in \FF$ such that $\theta^s \neq 1$, $1 \leq s \leq r$.	  
	$(E, D_{U_{k,n}}(i;\theta),f)$
	is a weighted $t$-design for $1 \leq i < n-k$, and $f(X)=p(U_{k,n}.X;\theta)$ for all $X \in D_{U_{k,n}}(i;\theta)$.
\end{Example}

\subsection{Further Implications}
We now obtain a weaker form of the Assmus-Mattson Theorem for matrix codes as a direct consequence of Theorem \ref{th:weightedtdesign}. Note that the result for subspace designs (those weighted designs for which $f(B)=1$ for every block $B$) obtained from rank-metric codes was shown in \cite{byrne2019assmus} with the further assumption that the number of codewords with a given support was dependent only on the dimension $i$ of that space for some range of $i$.

\begin{Corollary}\label{cor:AM}
	Let $C$ be an $\FF_{q}$-$[n\times m,k,d]$ rank-metric code. Let $t< d$ be a positive integer and let $C^\perp$
	have no more than $d-t$ distinct rank weights in the set $\{1,\ldots,n-t\}$.
	For each $i \in \{d,\ldots,n-t\}$, let $$B(i)=\{ U \leq E : \dim(U) = i, |C_{=U^\perp}| \neq 0 \}.$$
	Then for each $i \in \{d,\ldots,n-t\}$, $(E,B(i),f)$ is a weighted $t$-design over $\mathbb{F}_q$
	with $f(X):=|C_{=X^\perp}|$. 
\end{Corollary}

\begin{proof}
    Let $M:=M_C$.
	By Lemma \ref{lem:repqpoly}, we have that $M^* = M_{C^\perp}$ and for any $i \in \{1,\ldots,n\}$, $W_i(C^\perp) = A_{M^*}(i;q)$. 
	Also, $p(M.X;q) = |C_{=X^\perp}|$ for any subspace $X \leq E$. 
	Now $d_M = \min\{\dim X : X \text{ is a cocircuit of }M  \}$, which by Proposition \ref{prop:D_M(d_M;z)}, is the minimum dimension
	of any subspace $X$ such that $p(M.X;q) \neq 0$.  
	
	Since $C$ has minimum distance $d$, by Lemma \ref{lem:repqpoly} (4) there exists a $d$-dimensional subspace $X \leq E$ such that
	$|C_{=X^\perp}| = p(M.X;q) \neq 0$, while $p(M.U;q)=0$ for every subspace $U \leq E$ with $\dim(U) <d$. Therefore, $d=d_M$.
	By hypothesis, at most $d-t=d_M-t$ 
	of the integers $W_i(C^\perp)$ are non-zero for 
	$i \in \{1,\ldots ,n-t\}$. 
	By Lemma \ref{lem:repqpoly}, if $A_{M^*}(i;q) = 0 $, then $ A_{M^*/T}(i;q) = 0$,
	for any $t$-dimensional subspace $T \leq E$. 
	Therefore $M$ satisfies the hypothesis of 
	Corollary \ref{firstcorollary} and so the result follows. 	
\end{proof}

In the case of a $q$-matroid $M$ satisfying the hypothesis of Theorem \ref{th:weightedtdesign}, with an extra assumption on the cocircuits of $M$, our results imply the existence of a subspace design. These results form a direct $q$-analogue of the classical case (c.f. \cite[Section 3]{BRS2009}).

\begin{Lemma}\label{cor:pcocirc}
    Let $M$ be a $q$-matroid, $\theta \in \FF, \theta \neq 1$ and let $p$ be the greatest integer such that
    any subspace $X \leq E$  of dimension at most $p$ contains at most one cocircuit of $M$.
    Then, for each $i \in \{d_M,\ldots, p\}$, we have  
    $D_{M}(i;\theta) = \{ C \leq E : C \text{ is a cocircuit of }M, \dim(C)=i \}.$
\end{Lemma}

\begin{proof}
    If $C$ is a cocircuit of $M$, then $p(M.C;\theta)=\theta -1 \neq 0$ and so 
    $C \in D_{M}(\dim(C);\theta)$.
    Now let $X \in D_{M}(i;\theta)$ for some $i\leq p$. 
    Then $p(M.X;\theta) \neq 0$ and $X$ is a dependent space of $M^*$ of dimension at most $p$, so $X$ contains a unique circuit of $M^*$. By Theorem \ref{th:cocircs}, we have $X=C$ and the result follows. 
\end{proof}

\begin{Corollary}\label{cor:design}
    Let $M$ be a $q$-matroid that has at least one circuit and one cocircuit.
    Let $t<d_M$ be a positive integer such that the hypothesis of Theorem \ref{th:weightedtdesign} holds for some $\theta \in \FF,\theta \neq 1$. Let $p$ be the greatest integer such that
    any subspace $X \leq E$  of dimension at most $p$ contains at most one 
    cocircuit (respectively, at most one circuit) of $M$. Then for each $i \in \{d_M,\ldots,p\}$ (respectively, $\{d_{M^*},\ldots,p\}$) the set of cocircuits (respectively, the set of circuits) of $M$ of dimension $\min\{i,n-t\}$ forms the blocks of a $t$-subspace design. 
    Consequently, for each $i \in \{d_M,\ldots,p\}$ (respectively, $\{d_{M^*},\ldots,p\}$), the set of hyperplanes of $M$ (respectively, of $M^*$)
    of dimension $n-i$ forms the blocks of a $t$-subspace design.
\end{Corollary}

\begin{proof}
    From Corollary \ref{cor:pcocirc}, for each $i \in \{d_M,\ldots,p\}$ we have that $C_i:=D_M(i;\theta)$ is the set of cocircuits of $M$ of dimension $i$. Then by Corollary \ref{firstcorollary}, for each $i \in\{d_M,\ldots, p\}$, $C_i$ is the set of blocks of a weighted $t$-subspace design with $f(X) = p(M.X;\theta) = \theta-1$. Define a function $\hat{f}:D_M(i;\theta)\longrightarrow \FF$ by
    $\hat{f}(X) =(\theta-1)^{-1}f(X)$. This yields a $t$-subspace design ${\mathcal D}_i$ whose blocks are $C_i$. By \cite[Corollary 71]{bcj},
    for each $i$-dimensional cocircuit $X$ of $M$, $X^\perp$ is a hyperplane of $M$ and has dimension $n-i$.
    By Corollary \ref{cor:weighteddesign}, the set of hyperplanes of $M$ of dimension $n-i$ form the blocks of a $t$-subspace design, i.e., the dual design of ${\mathcal D}_i$.
    With the same arguments as above, by Corollary \ref{secondcorollary} the analogous statements hold for
    the circuits of $M$ and the hyperplanes of $M^*$. 
\end{proof}

An element $c$ of an $\FF_{q^m}$-$[n,k,d]$ vector rank-metric code $C$ is called {\em minimal} 
if for any $c' \in C$, $\sigma(c') \leq \sigma(c) \implies c' \in \langle c \rangle_{\FF_{q^m}}:=\{\nu c : \nu \in \FF_{q^m}\}$. In this case, for $U=\sigma(c)$, we have $p(M.X;q^m)=|C_{=U}| = q^m-1$. If every codeword of rank $i$ in $C$ is minimal, then $A_{M_{C^\perp}}(i;q^m) = W_i(C^\perp) = (q^m-1)|D_M(i;q^m)|$.
If we apply this with Corollary \ref{cor:design}, we retrieve the Assmus-Mattson Theorem for 
$\FF_{q^m}$-$[n,k,d]$ codes (c.f. \cite{byrne2019assmus}).

\begin{Corollary}\label{cor:codetodesign}
    Let $C$ be an $\FF_{q^m}$-$[n,k,d]$ code. Let $t< d$ be a positive integer and let $C^\perp$
    be an $\FF_{q^m}$-$[n,n-k,d^\perp]$ code
	having no more than $d-t$ distinct rank weights in the set $\{1,\ldots ,n-t\}$.
	Let $p$ be the greatest integer such that
	every codeword of $C$ of rank at most $p$ is minimal.

	\begin{enumerate}
	    \item The supports of the words of rank weight $d$ in $C$ (respectively $d^\perp$ in $C^\perp$) form the blocks of a $t$-design over $\FF_q$.
	    \item For each $i \in \{d,\ldots,p\}$ (respectively, $\{d^\perp,\ldots,p\}$) the supports of the minimal codewords of $C$ (respectively $C^\perp)$ of dimension $\min\{i,n-t\}$ form the blocks of a $t$-design over $\FF_q$.
	\end{enumerate}
\end{Corollary}

\begin{Example}
    In \cite[Theorem 12]{randrianarisoa}, it is shown that any non-degenerate $\FF_{q^m}$-$[N,k>1]$ rank-metric code with constant weight $d$ satisfies $N=km$, $d=m$ and is generated by a matrix 
    $G \in \FF_{q^m}^{k \times N}$ whose $N$ columns form a basis of $\FF_{q^m}^k$ as an $\FF_q$-vector space.
    Moreover, the dual code has minimum distance $2$. 
    Let $C^\perp$ be an $\FF_{q^m}$-$[km,k,m]$ constant weight code constructed as above. Let $M=M_C$, so
    that $M^*=M_{C^\perp}$. For any $X \leq \FF_q^{km}$, we have $p(M.X;q^m) = 0$ unless $X$ is the support of a codeword
    of $C^\perp$, in which case $\dim(X) = m$. Therefore, 
    $A_{M^*}(m;q^m) = q^{km}-1, A_{M^*}(0;q^m) = 1$ and $A_{M^*}(i;q^m)=0$ for $i \neq 0,m$.
    Then $d_{M} =d= 2$ and $R_{M}(2;\theta)=\{ m\}$. 
    Therefore, by Corollary \ref{cor:codetodesign} the cocircuits of $M$ of dimension 2, which are the supports of codewords of rank $2$, form a $1$-design over $\FF_q$.
    Similarly, the supports of the words of rank $m$ in $C^\perp$ form the blocks of a 1-design over $\FF_q$, in fact a $1$-$(km,m,1;q)$ design, which is a $q$-Steiner system, whose blocks form a spread in $\FF_{q}^{km}$.
\end{Example}

While Theorem \ref{th:weightedtdesign} has considerable potential for constructing weighted subspace designs, utilizing it requires constructions of a $q$-polymatroid $M$ whose weight enumerator takes few non-zero values and whose cocircuits have large enough dimension. Most $q$-matroids and $q$-polymatroids are not representable, however those that are, i.e. those that can be represented by rank metric codes, offer more tangible constructions.

 In order to search for examples of rank-metric codes satisfying the conditions in Corollary \ref{cor:codetodesign} we implemented in Magma~\cite{Magma} a search through 
  random $\mathbb{F}_{2^m}$-$[n,k,d]$ rank metric codes, for different values of $m$, $n$, and $k$ and $t=2$. 
  We make some remarks on the parameters of potentially interesting codes. 
  
  A matrix code is called maximum rank distance (MRD) if it has parameters $\Fq$-$[n \times m,k,d]$ with $k=\max\{m,n\}(\min\{n,m\}-d+1)$. The MRD $\Fqm$-linear codes have parameters
  $\Fqm$-$[n,n-d+1,k]$. MRD codes do satisfy the criteria of Corollaries \ref{cor:AM} and \ref{cor:codetodesign} and there are several constructions of them. However, the corresponding subspace designs associated with these codes are trivial. We therefore would exclude them from our search space.
  Note that, as either $m$ or $q$ grow asymptotically, the $\mathbb{F}_{q^m}$-linear MRD codes are dense in the space of linear codes (see for example \cite{ByrneRavagnani,Neri+}), just as MDS codes are dense as $q$ becomes large. However, for small values of $q$ and  $m$,  we do not necessarily have a high probability of a random code being MRD, and this was checked and confirmed for all sets of parameters chosen for experiments. 
  We also have modest constraints on the sizes of $m,n$ to exclude trivialities.
  For example with $t=2$, we require that $d \geq 3$. If $m \leq 4$ and $d=3$, then to find a suitable code $C$, we would require $C^{\perp}$ to be a one-weight code.
  However, by \cite[Proposition 4.6]{byrne2019assmus}, 
  the dual code of a constant weight $\Fq$-$[n \times m,k]$ code has minimum distance at most $2$.
  If $m=4=d$, then $C$ must be a one-weight code and so $C^\perp$ has weight at most $2$, which would not yield interesting results.
  We therefore set $m>4$
  in order to meet the criteria of Corollary \ref{cor:codetodesign}. This means that we search for linear codes over alphabets of size at least $2^5$. In most experiments, we chose $m=n-1$ or $m=n$ to increase the probability of satisfying the criteria.

 Each code is  given by a generator matrix in standard form for a linear code, i.e.  $(I_k | A)$, where $A$  goes through the space of $k \times (n-k)$-matrices with entries in $\mathbb{F}_2^m$, up to equivalence under the action of the Galois group Gal$(\mathbb{F}_{2^m}/\mathbb{F}_2)$. This yields a search space of size $2^{m(k(n-k)-1)}$ matrices, which is quickly out of reach of a computer, even for small values of $k$ and $n$. Ideally, a single representative in each equivalence class for the underlying $q$-polymatroids should be computed, but it is not clear to us how to pre-compute these representatives such that running the search would be more time effective.   
 
In our algorithm, we first compute the weight distribution of the rank metric code by going through all code words (up to a scalar) and then we deduce the weight distribution of the dual code by using the MacWilliams identities. The code is publicly available ~\cite{code}. For each set of parameters, we ran the code on a different core of an 40GHz Intel Xeon E5-2640 processor, and we set a timeout of 16 days for each run.
The number of codes that we were able to check in this way are given in the fifth column of Table~\ref{experiments}.  

\begin{table}[h!]
    \centering
    \begin{tabular}{|c|c|c|c|c|}
      \hline 
      $m$ & $n$ & $k$  & no. of codes checked & \% of search space \\
      \hline 
      5  &  6 & 2 &  405,285,656 &  $0.125$\\
      \hline 
      6 & 6 & 2 &   146,666,189 & $3.334 \times 10^{-5}$\\
      \hline
      6 & 6 & 3 &  442,349 &  $1.572 \times 10^{-9}$  \\
      \hline
      6 & 8 & 2 &  44,700,000 & $6.058 \times 10^{-13}$\\
      \hline
      7 & 8 & 2 & 13,800,000  & $9.132 \times 10^{-17}$\\
      \hline
      8 & 8 & 2 &   3,800,000 & $1.228 \times 10^{-20}$ \\
       \hline
    \end{tabular}
    \caption{Random search through $\mathbb{F}_{q^m}$-$[n,k,d]$ random rank metric codes}
    \label{experiments}
\end{table}

 While a complete search is far from complete, these numbers suggest that a more systematic search for higher parameter values would be needed to effectively construct examples of rank metric codes yielding weighted subspace designs. More generally, what is really required is a theoretical approach to construct rank-metric codes and $q$-(poly)matroids with prescribed weight distributions.  

\vspace{0.5 cm}

\section{Acknowledgements}

The authors would like to say thanks to Elif Sa\c{c}ikara for many useful discussions in the preparation of this work. The authors are grateful to Heide Gluesing-Luerssen and Benjamin Jany for their helpful remarks on $q$-polymatroid equivalence.
This paper stems from a collaboration that was initiated at the Women in Numbers Europe (WIN-E3) conference, held in Rennes, August 26-30, 2019. The authors are very grateful to the organisers: Sorina Ionica, Holly Krieger, and Elisa Lorenzo Garc\'ia, for facilitating their participation at this workshop, which was supported by the Henri Lebesgue Center, the Association for Women in Mathematics (AWM) and the Clay Mathematics Institute (CMI).

\bibliographystyle{abbrv}
\bibliography{References}

\begin{thebibliography}{10}

\bibitem{AM69}
E.~F. Assmus, Jr. and H.~F. Mattson, Jr.
\newblock New {$5$}-designs.
\newblock {\em J. Combinatorial Theory}, 6:122--151, 1969.

\bibitem{Birkhoff}
G.~Birkhoff.
\newblock {\em Lattice Theory}.
\newblock American Mathematical Society, revised edition, 1948.

\bibitem{Magma}
W.~Bosma, J.~Cannon, and C.~Playoust.
\newblock The {M}agma algebra system. {I}. {T}he user language.
\newblock {\em J. Symbolic Comput.}, 24(3-4):235--265, 1997.
\newblock Computational algebra and number theory (London, 1993).

\bibitem{bks18}
M.~Braun, M.~Kiermaier, and A.~Wassermann.
\newblock $q$-{A}nalogues of designs: Subspace designs.
\newblock In M.~Greferath, M.~O. Pav{\v{c}}evi{\'c}, N.~Silberstein, and
  M.~{\'A}. V{\'a}zquez-Castro, editors, {\em Network Coding and Subspace
  Designs}, pages 171--211. Springer, 2018.

\bibitem{BRS2009}
T.~Britz, G.~Royle, and K.~Shiromoto.
\newblock Designs from matroids.
\newblock {\em SIAM J. Discrete Math.}, 23(2):1082--1099, 2009.

\bibitem{britzshiromotomacwill}
T.~Britz and K.~Shiromoto.
\newblock A {M}ac{W}illiams type identity for matroids.
\newblock {\em Discrete Math.}, 308(20):4551--4559, 2008.

\bibitem{code}
E.~Byrne, M.~Ceria, S.~Ionica, and R.~Jurrius.
\newblock Weighted designs from rank metric codes.
\newblock \url{https://gitlab.com/ISorina/weighted-designs}.

\bibitem{WINEpaper1}
E.~Byrne, M.~Ceria, S.~Ionica, R.~Jurrius, and E.~Sa\c{c}{\i}kara.
\newblock Constructions of new matroids and designs over $\mathbb{ F}_q$.
\newblock {\em Designs, Codes and Cryptography}, 2022.

\bibitem{bcj}
E.~Byrne, M.~Ceria, and R.~Jurrius.
\newblock Constructions of new {$q$}-cryptomorphisms.
\newblock {\em J. Combin. Theory Ser. B}, 153:149--194, 2022.

\bibitem{byrne2019assmus}
E.~Byrne and A.~Ravagnani.
\newblock An {A}ssmus-{M}attson theorem for rank metric codes.
\newblock {\em SIAM J. Discrete Math.}, 33(3):1242--1260, 2019.

\bibitem{ByrneRavagnani}
E.~Byrne and A.~Ravagnani.
\newblock Partition-balanced families of codes and asymptotic enumeration in
  coding theory.
\newblock {\em {Journal of Combinatorial Theory} Series A}, 171, 2020.

\bibitem{Cameron}
P.~J. Cameron.
\newblock Generalisation of {F}isher's inequality to fields with more than one
  element.
\newblock In {\em Combinatorics ({P}roc. {B}ritish {C}ombinatorial {C}onf.,
  {U}niv. {C}oll. {W}ales, {A}berystwyth, 1973)}, London Math. Soc. Lecture
  Note Ser., No. 13, pages 9--13. Cambridge Univ. Press, London, 1974.

\bibitem{qdist}
C.~A. Charalambides.
\newblock {\em Discrete {$q$}-distributions}.
\newblock John Wiley \& Sons, Inc., Hoboken, NJ, 2016.

\bibitem{CrapoRota}
H.~H. Crapo and G.-C. Rota.
\newblock On the foundations of combinatorial theory. {II}. {C}ombinatorial
  geometries.
\newblock {\em Studies in Appl. Math.}, 49:109--133, 1970.

\bibitem{del76}
P.~Delsarte.
\newblock Association schemes and {$t$}-designs in regular semilattices.
\newblock {\em J. Combinatorial Theory Ser. A}, 20(2):230--243, 1976.

\bibitem{DL17}
C.~Ding and C.~Li.
\newblock Infinite families of 2-designs and 3-designs from linear codes.
\newblock {\em Discrete Math.}, 340(10):2415--2431, 2017.

\bibitem{GLJindep}
H.~Gluesing-Luerssen and B.~Jany.
\newblock Independent spaces of $q$-polymatroids.
\newblock \url{https://arxiv.org/abs/2105.01802}, 2021.

\bibitem{GLJ}
H.~Gluesing-Luerssen and B.~Jany.
\newblock $q$-{P}olymatroids and their relation to rank-metric codes.
\newblock {\em Journal of Algebraic Combinatorics}, pages 1--29, 2022.

\bibitem{gorla2019rank}
E.~Gorla, R.~Jurrius, H.~H. L\'{o}pez, and A.~Ravagnani.
\newblock Rank-metric codes and {$q$}-polymatroids.
\newblock {\em J. Algebraic Combin.}, 52(1):1--19, 2020.

\bibitem{JP18}
R.~Jurrius and R.~Pellikaan.
\newblock Defining the {$q$}-analogue of a matroid.
\newblock {\em Electron. J. Combin.}, 25(3):Paper No. 3.2, 32, 2018.

\bibitem{KP}
M.~Kiermaier and M.~O. Pav\v{c}evi\'{c}.
\newblock Intersection numbers for subspace designs.
\newblock {\em J. Combin. Des.}, 23(11):463--480, 2015.

\bibitem{kung1996critical}
J.~P.~S. Kung.
\newblock Critical problems.
\newblock In {\em Matroid theory ({S}eattle, {WA}, 1995)}, volume 197 of {\em
  Contemp. Math.}, pages 1--127. Amer. Math. Soc., Providence, RI, 1996.

\bibitem{zhizheng}
M.~Liu and Z.~Zhang.
\newblock {$q$}-analog of determinant of a kind of binomial coefficient
  matrices.
\newblock {\em Math. Morav.}, 8(1):15--24, 2004.

\bibitem{MT18}
J.~V.~S. Morales and H.~Tanaka.
\newblock An {A}ssmus-{M}attson theorem for codes over commutative association
  schemes.
\newblock {\em Des. Codes Cryptogr.}, 86(5):1039--1062, 2018.

\bibitem{Neri+}
A.~Neri, A.-L. Horlemann-Trautmann, T.~Randrianarisoa, and J.~Rosenthal.
\newblock On the genericity of maximum rank distance and gabidulin codes.
\newblock {\em {Des. Codes Cryptogr. }}, 86:341–363, 2018.

\bibitem{oxley}
J.~Oxley.
\newblock {\em Matroid theory}, volume~21 of {\em Oxford Graduate Texts in
  Mathematics}.
\newblock Oxford University Press, Oxford, second edition, 2011.

\bibitem{randrianarisoa}
T.~H. Randrianarisoa.
\newblock A geometric approach to rank metric codes and a classification of
  constant weight codes.
\newblock {\em Des. Codes Cryptogr.}, 88(7):1331--1348, 2020.

\bibitem{SKH04}
D.-J. Shin, P.~V. Kumar, and T.~Helleseth.
\newblock An {A}ssmus-{M}attson-type approach for identifying 3-designs from
  linear codes over {$Z_4$}.
\newblock {\em Des. Codes Cryptogr.}, 31(1):75--92, 2004.

\bibitem{shiromoto19}
K.~Shiromoto.
\newblock Codes with the rank metric and matroids.
\newblock {\em Des. Codes Cryptogr.}, 87(8):1765--1776, 2019.

\bibitem{suzuki1990inequalities}
H.~Suzuki.
\newblock On the inequalities of {$t$}-designs over a finite field.
\newblock {\em European J. Combin.}, 11(6):601--607, 1990.

\bibitem{Lint}
J.~van Lint and R.~Wilson.
\newblock {\em A Course in Combinatorics}.
\newblock Cambridge {University Press}, 1992.

\bibitem{welsh}
D.~J.~A. Welsh.
\newblock {\em Matroid theory}.
\newblock L. M. S. Monographs, No. 8. Academic Press [Harcourt Brace
  Jovanovich, Publishers], London-New York, 1976.

\bibitem{whittle}
G.~Whittle.
\newblock Characteristic polynomials of weighted lattices.
\newblock {\em Adv. Math.}, 99(2):125--151, 1993.

\bibitem{zaslavsky1987mobius}
T.~Zaslavsky.
\newblock The {M}{\"o}bius function and the characteristic polynomial.
\newblock {\em Combinatorial geometries}, 29:114--38, 1987.

\end{thebibliography}

\end{document}